\documentclass[11pt]{amsart}
\usepackage{amssymb, mathtools, latexsym, mathrsfs,   color }
\usepackage[all]{xy}
\usepackage[colorlinks=true, pdfstartview=FitV,linkcolor=blue,citecolor=blue,urlcolor=blue]{hyperref}

    \advance \topmargin by -\headheight
    \advance \topmargin by -\headsep

    \evensidemargin \oddsidemargin
 \newtheorem{theorem}{Theorem}[section]
  \newtheorem{prop}[theorem]{Proposition}
  
  \newtheorem{lemma}[theorem]{Lemma}

\newcommand{\la}{\langle}
\newcommand{\ra}{\rangle}

\newcommand{\Sp}{{\rm Sp}}

\newcommand{\deff}{\stackrel{\rm def}{=}}

\newcommand{\ad}{{\rm ad}}
\newcommand{\St}{\bold{St}}
\newcommand{\ol}{\overline}

\renewcommand{\rm}{\mathrm}

\newcommand{\cal}{\mathcal}

\theoremstyle{definition}

\newcommand{\matrixx}[4]{\begin{pmatrix}
#1 & #2 \\ #3 & #4
\end{pmatrix} }

\numberwithin{equation}{section}

\begin{document}

\title{Weil Representations of  Twisted Loop Groups of Type $A_n^{(2)}$ }

\author[Yanze Chen,  Yongchang Zhu]{Yanze Chen$^*$ and Yongchang Zhu}

\address{Department of Mathematics, The Hong Kong University of Science and Technology, Clear Water Bay, Kowloon, Hong Kong}

\email{ychenen@connect.ust.hk}
\email{mazhu@ust.hk}

\thanks{$*$This research is supported by Hong Kong RGC grant 16301718.}


\begin{abstract}
We construct Weil representations of twisted loop groups  of type $A_n^{(2)}$ over local fields. We prove that the associated cover of the twisted loop group
  is the two fold metaplectic cover of the affine Kac-Moody group of type $A_n^{(2)}$ given by Patnaik-Puskas.
\end{abstract}

\maketitle

\section{Introduction}

The Weil representation of symplectic groups plays an important role in the explicit constructions of automorphic forms.
A generalization of Weil representation of symplectic loop groups  over local fields, which has Lie algebra of type $C_n^{(1)}$,  was given in \cite{Z}.  As in the classical case,
it has applications to the theory of  automrphic forms on loop groups  \cite{GZ1} \cite{GZ2} \cite{LZ}. A less known fact is that
  a loop group of type $A_n^{(2)}$ has a realization as a symplectic group. The purpose of this work is to construct
   its Weil representation and prove that the associated cover is  the two fold metaplectic cover of the affine Kac-Moody group of type $A_n^{(2)}$ given by Patnaik-Puskas \cite{PP}.

To give a realization of  twisted group of type $A_{n-1}^{(2)}$ over a field $F$ as a symplectic group,
 we introduce  a symplectic form on the $n$-dimensional space  $ F^{n} (( t)) $ over the Laurent power series field $F(( t )) $
 by
\begin{equation}\label{1}
  \la  u ( t ) , v ( t ) \ra  =   {\rm Res} (  u (  t )  C v ( - t )^T ) .
  \end{equation}
 where $C$ is an $ n \times n$ symmetric non-degenerate matrix over $F$, and ${\rm Res} (f)$ for  $ f\in F(( t ))$ is equal to the coefficient of $t^{-1}$ in $f$.  The form is $F$-valued.
 We assume that the vectors in $ F^{n} ((t))$ are row vectors and the group  $GL_{n} ( F (( t)) )$
 acts on $ F^{n} ((t))$ by the right multiplication.
 We are interested in the group that consists of the $F((t))$-linear isomorphisms of $ F(( t))^n $ that preserves the symplectic form
  $\la \; , \; \ra $.
A $F((t))$-linear isomorphism  $ g ( t ) \in GL_n ( F (( t)) )$ preserves $\la \; , \; \ra $ iff
\begin{equation}
  g(  t )  C  g ( - t )^T  = C
\end{equation}
which is equivalent to
\[         C  g ( -t )^{-T} C^{-1}  =  g ( t ) \]
i.e.,  $g(t)$ is in the fixed point of the automorphism
\[ \sigma :  g ( t) \mapsto   C g ( -t )^{-T} C^{-1}  .\]
Since $C$ is symmetric, we have $ C^{-1} C^T = C^{-1} C = I $, so  $ \sigma $ has order $2$.
A $F((t))$-isomorphism of $F((t))^n$ preserving symplectic form (1.1)
   is an element $  g \in GL_n ( F(( t)))$ satisfying $ \sigma ( g ) = g$.
  The twisted loop group $LG(F)$ is
\begin{equation}
   LG(F) = \{   g \in SL_n ( F(( t))) \, | \, \sigma ( g ) = g \}.
 \end{equation}
 A central extension of $LG(F)$ is the Kac-Moody group of type $A_{n-1}^{(2)}$ (see Section 2).

 The formal Weil representation of twisted groups $LG({\Bbb C} )$ is given in \cite{Gi} and applied to the study of Gromov-Witten invariants.
The Weil representation studied in this work is of analytic nature and is similar to the approach in \cite{Z}.
We use Lagrangian subspace decomposition
\begin{equation}
  F(( t))^n = F [[ t]]^n  \oplus F [ t^{-1} ]^n t^{-1}
\end{equation}
The general result \cite{Z} implies that there is a projective representation of $LG(F)$ on
the space ${\cal S} (  F [ t^{-1} ]^n t^{-1} )$ of schwartz functions on $F [ t^{-1} ]^n t^{-1} $.
However, our formulation here is different from \cite{Z} that the symbols associated to
 the projective representation are treated differently.
In \cite{Z}, the loop symplectic group $Sp_{2n}(F((t)))$ is viewed as the Chevalley group of
type $C_n$ over the field $F((t))$, the symbol
 is computed as a symbol over the field $F((t))$. This approach seems not work for twisted loop
  group $A_n^{(2)}$, as it is not a split group over $F((t))$.
  We will adapt Patnaik-Puskas's approach to metapectic covers of Kac-Moody groups, which generalize the work of Mutsumoto \cite{Ma}  about metaplectic cover for finite dimensional semisimple groups to Kac-Moody groups.
 Our main result can be described as a construction of Weil representation of metaplectic cover of Kac-Moody group
  over local fields associated to the Hilbert symbol.

   This paper is organized as follows. In Section 2, we will give a realization the affine Kac-Moody algebra of type $A_{n-1}^{(2)}$
    as a symplectic Lie algebra, corresponding to the realization of twisted group in (1.3). We fix a set of
     Chevalley generators of the Lie algebra and fix an explicit root datum used later and a Weyl group invariant quadratic form in the coroot lattice.
   In Section 3, we will construct the Weil representation of
the twisted loop group as projective representation following \cite{Z}.
In Section 4, we recall Tits' Kac-Moody group functor and the construction of metaplectic cover of Kac-Moody groups by Patnaik-Puskas  .
 We give a description of the Patnaik-Puskas's cover in terms of generators and relations.
 In Section 5, we show that the Weil representation in Section 3 is in fact a representation of
  the metaplectic cover of Kac-Moody group over a local field of type $A_n^{(2)}$ associated to the Hilbert symbol.
  Our proofs use the generators and relations obtained in Section 4.

\

 \section{\Large Symplectic Realization of Affine Kac-Moody algebras of type $A_{n}^{(2)}$.   }

In this section we give a symplectic realization of affine algebra of type $A_n^{(2)}$ corresponding to the twisted loop
group constructed in Section 1. We then fix a simply connected datum root for $A_n^{(2)}$ using the realization.

The affine algebra of type $A_{n}^{(2)}$ will be constructed as a subalgebra of
\[ \widehat{sl}_{n+1} = sl_{n+1} [ t , t^{-1} ] \oplus {\Bbb C} K \oplus {\Bbb C} d \]
 with Lie bracket
\[    [ a t^m , b t^l ] = [ a , b ] t^{m+l} + m \delta_{ m + l , 0 }  {\rm Tr} ( a b ) K ,   \; \; [ d , a t^m ] = m a t^m .\]
It has a symmetric invariant bilinear form given by
\begin{equation}\label{2.1}
 ( a t^m , b t^n ) = {\rm Tr} ( ab ) \delta_{m, - n } ,  (a t^m , K ) = ( at^m , d ) = ( K , K ) = ( d , d ) = 0 , ( K , d ) = 1.
\end{equation}

Let $C$ be an $(n+1)\times (n+1)$ symmetric non-generate matrix with integer entries.
We use $C$ to define a non-degenerate symplectic form on ${\Bbb C} [ t , t^{-1} ]^{n+1} $ as in (1.1).
The Lie algebra $sl_{n+1} [ t , t^{-1} ]$ acts on ${\Bbb C} [ t , t^{-1} ]^{n+1} $ from the right.
   It is easy to see that $ a( t) \in sl_{n+1}{\Bbb C} [ t , t^{-1} ]^{n+1} $ preserves (1.1) iff
    $  a ( t ) =   C a( - t)^{-T} C^{-1} $, such elements form a subalgebra of $sl_{n+1} [ t , t^{-1} ]$.
    This subalgebra gives arise a subalgebra of $\widehat{sl}_{n+1}$ of type $A_n^{(2)}$. To be more precise, we first define an automorphism of ${\frak g} = sl_{n+1} $ by
\begin{equation}
\sigma ( a ) = - C a^T C^{-1} ,
\end{equation}
 The fixed point algebra
  \[   {\frak g}_0 = \{ a \in {\frak g} \, | \, \sigma ( a) = a \}  , \]
 is the orthogonal Lie algebra preserving the symmetric bilinear form
 \begin{equation}
   xC x^T .
 \end{equation}
 Then we extend the automorphism $\sigma $ to  $sl_{n+1} [ t , t^{-1} ]$ by
 \begin{equation}    \sigma ( at^n ) =  \sigma ( a) ( -t)^n \end{equation}
It is clear that  $\sigma $ commutes with $ d $.
It is easy to check that  $\sigma $ preserves  the cocycle
 \[ ( a ( t) , b ( t ) ) =  {\rm Res}  \left(     {\rm tr} (  a' ( t) b ( t ) ) \right) . \]
So $\sigma $ extends to an order $2$ automorphism on $\widehat{sl}_{n+1}$ by
\[ \sigma ( d ) = d , \; \; \; \sigma ( K ) = K  .\]
Let
\[ {\frak g}_1 = \{ a \in {\frak g} \, | \, \sigma ( a) = -a \} ,\]
 which is an irreducible module over ${\frak g}_0$.
 The fixed point Lie algebra $\widehat{sl}_{n+1}^{\sigma } $ is
\[ \widehat{sl}_{n+1}^{\sigma } = {\frak g}_0 [ t^2 , t^{-2} ] + {\frak g}_1 [ t^2 , t^{-2} ] t +  {\Bbb C} K +  {\Bbb C} d .\]
Our explicit description of a set of simple roots and coroots  will prove that  it has type
 $ A_{n}^{ (2 )} $.   Note our ordering of simple roots for $A_{2l}^{(2)} $ given in the later part of this section is the reverse
   of that in \cite{K}.
 The realization of given in \cite{K} uses a different order $2$ automorphism $\widehat{sl}_{n+1}^{\sigma } $, where the affine Lie algebra is realized as a central extension of
  an orthogonal Lie algebra.

  Since we want our affine Kac-Moody group to be split,
 we take $C$ as follows. For $sl_{2l+1}$, $ l \geq 1 $, we take
\[   C = \left( \begin{matrix}  0 & I_l & 0  \\   I_l & 0 & 0 \\ 0 & 0 & 1  \end{matrix} \right) .\]
The fixed point algebra will be $ A_{2l}^{ ( 2)} $.

  For $sl_{2l}$, $ l \geq 1 $,
    \[   C = \left( \begin{matrix}  0 & I_l \\   I_l & 0 \end{matrix} \right) .\]
 The fixed point algebra will be $ A_{2l-1 }^{ ( 2)} $.
  Notice that in both cases, we have $C= C^T = C^{-1} $.

\

\noindent {\bf 1. Symplectic Realization of $A_{2}^{(2)}$.}
We have ${\frak g}= sl_3 $. The fixed point Lie algebra ${\frak g}_0$ is the orthogonal Lie algebra preserving the bilinear form (2.2)
 on ${\Bbb C}^3 $, so
it is isomorphic to $sl_2 $.
More explicitly, we have
 \[ {\frak g}_0 = {\Bbb C} h_1 +{\Bbb C} e_1 + {\Bbb C} f_1  , \]
where
   \[  h_1 = \left( \begin{matrix}  2 & 0 & 0  \\   0 & -2 & 0 \\ 0 & 0 & 0  \end{matrix} \right),
 \; e_1 = \left( \begin{matrix}  0 & 0 & 2  \\   0 & 0 & 0 \\ 0 & -2 & 0  \end{matrix} \right) , \;
   f_1 = \left( \begin{matrix}  0 & 0 & 0  \\   0 & 0 & -1 \\ 1 & 0 & 0  \end{matrix} \right) .\]
  They form a standard $sl_2$-triple. We write $h_1 = 2 \epsilon_1 $, $\eta_0 = {\Bbb C} \epsilon_1 $.
  We identify the dual space $\eta_0^*$ with $\eta_0$ by the bilinear form $ ( \epsilon_1 , \epsilon_1 ) =1$.
   The $(-1)$-eigenspace of $\sigma $ is
 \[  {\frak g}_1 = \rm{Span} ( v =\left( \begin{matrix}  0 & 1 & 0  \\   0 & 0 & 0 \\ 0 & 0 & 0  \end{matrix} \right),
  \;  \left( \begin{matrix}  0 & 0 & 1  \\   0 & 0 & 0 \\ 0 & 1 & 0  \end{matrix} \right), \;
    \left( \begin{matrix}  1 & 0 & 0  \\   0 & 1 & 0 \\ 0 & 0 & -2  \end{matrix} \right), \;
   \left( \begin{matrix}  0 & 0 & 0  \\   0 & 0 & 1 \\ 1 & 0 & 0  \end{matrix} \right), \;
  u=  \left( \begin{matrix}  0 & 0 & 0  \\   1 & 0 & 0 \\ 0 & 0 & 0  \end{matrix} \right) )  \]
 They have eigenvalues $ 4, 2, 0 , -2 , -4 $ under $ {\rm ad} \, h_1 $.
 \[   e_0= u t,  \; f_0= v t^{-1} ,   h_0 = K - \frac 1 2 h_1   \]
form a standard triple of $sl_2 $.
The Chevalley generators of Lie algebra $A_2^{(2)}$  are
\begin{equation}  e_0 , f_0 , h_0 , e_1 , f_1 , h_1 \end{equation}
and we fix a root datum $\{ \Lambda , \Lambda^{\vee} , \Pi , \Pi^\vee \} $ as follows.
\[ \Pi^{\vee } = \{  \alpha_0^\vee,    \alpha_1^\vee \} \]
where $\alpha_0^\vee = K - \epsilon_1 $, $\alpha_1 ^\vee = h_1 = 2 \epsilon_1 $, $\Pi $ is the set of simple roots.
\[ \Lambda^\vee = {\Bbb Z}\alpha_0^\vee + {\Bbb Z} \alpha_1^\vee + {\Bbb Z} ( 2 d ) \]
Let $ \Lambda$ be the dual lattice sitting in the dual of Cartan subalgebra. 
Let $W$ be the Weyl group of the fixed point Lie algebra $\widehat{sl}_{3}^\sigma$. 
We define a $W$-invariant form
\[  Q :  \Lambda^{\vee} \times  \Lambda^{\vee} \to {\Bbb Z} \]
by
\[  Q ( \alpha_0^\vee , \alpha_0^\vee ) =1 , Q ( \alpha_1^\vee , \alpha_1^\vee ) = 4,   Q ( \alpha_0^\vee , \alpha_1^\vee ) = -2 ,
  Q ( 2d ,  \alpha_0^\vee ) =2 ,   Q ( 2d ,  \alpha_1^\vee ) =Q ( 2d , 2d ) =0 .\]
This form is $\frac 12 $ times the restriction of (2.1) on $\Lambda^\vee $.

\

\noindent {\bf 2. Symplectic Realization of $A_{2l}^{(2)}$ ($l\geq 2$).}

We have ${\frak g}  = sl_{2l+ 1} $ in this case. The $\sigma $-fixed Lie algebra  $ {\frak g}_0  $
is the orthogonal Lie algebra of type $B_l $, it has  Cartan subalgebra
\begin{equation}
   \eta_0 = \{ \rm{diag} ( d_1 , \dots , d_l , - d_1 , \dots , -d_l ) \}.
\end{equation}
We write
\[ \rm{diag} ( d_1 , \dots , d_l , - d_1 , \dots , -d_l )= d_1 \epsilon_1  + \dots + d_l \epsilon_l ,\]
 so  $ \epsilon_i $ is a diagonal matrix with $i$-th diagonal $1$, $(l+i)$-th diagonal $-1$ and all other diagonals $0$.
 The elements $\epsilon_1 , \dots , \epsilon_l$ form a basis for $\eta_0$. We introduce a symmetric bilinear form
  on $\eta_0$ by
  \[ ( \epsilon_i , \epsilon_j ) = \delta_{ij} ,\] Let $\epsilon_1^*,\cdots,\epsilon_n^*$ be the dual basis of $\epsilon_1,\cdots,\epsilon_n$ in $\eta_0^*$, the dual space of $\eta_0$. 
The roots of ${\frak g}_0 $ are
\[ \Delta_0 \deff \{  \pm \epsilon_i \pm \epsilon_j ,    1\leq  i \ne j \leq l \} \sqcup \{ \pm \epsilon_i ,  1\leq  i \leq l\}.  \]
 A set of simple roots is
 \[    \Pi_0 = \{ \alpha_{1}= \epsilon_1 - \epsilon_2, \alpha_{2} = \epsilon_{2}-\epsilon_3, \dots ,
 \alpha_{l-1}=\epsilon_{l-1} - \epsilon_{l} , \alpha_l = \epsilon_l \} \]
${\frak g}_1$ is
\[  {\frak g}_1 = \{  \left( \begin{matrix}  a_1 & a_2  & a_3  \\   b_1  & a_1^T & b_3
 \\ b_3^T & a_3^T  & c_3  \end{matrix} \right) \; | \;  a_2^T= a_2 , b_1^T= b_1, \; c_3= -2 {\rm tr} ( a_1)  \} \]
 As a $ {\frak g}_0$-module, it has weights
\[  \Delta_1 \deff  \{  \epsilon_i - \epsilon_j ,    1\leq  i \ne j \leq l \}
     \sqcup  \{ \pm( \epsilon_i + \epsilon_j ),    1\leq  i ,j \leq l \}   \sqcup \{ \pm \epsilon_i ,  1\leq  i \leq l\}.           \]
  ${\frak g}_1$ is a highest weight module of $ {\frak g}_0$-module with highest weight $ 2 \epsilon_1 $.

\

 The fixed point algebra $ \widehat {sl}_{2l+1}^\sigma $ is
 \begin{equation}
    \widehat {sl}_{2l+1}^\sigma = {\frak g}_0 \otimes {\Bbb C} [ t^2 , t^{-2} ] + {\frak g}_1 \otimes {\Bbb C} [ t^2 , t^{-2} ] t
      + {\Bbb C} K + {\Bbb C} d
 \end{equation}
 with Cartan subalgebra
 \begin{equation} \eta =  \eta_0 + {\Bbb C} K + {\Bbb C} d  .\end{equation}
The real roots of  $A_{2l}^{(2)}$ are
  \[ \Phi_{\rm re} = \{   \alpha + 2 k \delta \; | \;  \alpha \in \Delta_0, k \in {\Bbb Z} \}
   \sqcup \{   \alpha + ( 2 k -1) \delta \; | \;  \alpha \in \Delta_1, k \in {\Bbb Z} \}\]
   where $\delta\in\eta^*$ is the element with $\delta(\eta_0)=\delta(K)=0,\delta(d)=1$. 
   The set of simple roots is
   \[\Pi=\{\alpha_0=\delta-2\epsilon_1,\alpha_1, \dots , \alpha_l \}.\]
The corresponding set of simple coroots is
\[ \Pi^\vee=\{\alpha_0^\vee=K-\epsilon_1, \alpha_1^\vee = \epsilon_1 -\epsilon_2 ,
  \dots , \alpha_{l-1}^\vee = \epsilon_{l-1} -\epsilon_l ,  \alpha_l^\vee = 2 \epsilon_l \} \]

\

 We fix Chevalley generators. For $1\leq i \leq l-1$,
 \[    e_i = E_{ i ,  i +1  } -  E_{ l + i + 1 , l + i    }  , \; \;   f_i =  E_{ i +1 , i  } -  E_{ l + i  ,  l + i +1    } .\]
 For $ \alpha_l = \epsilon_l $,
 \[  e_l = 2 E_{l , 2 l+1 } - 2 E_{ 2 l+1 , 2 l }  , \; \; \;  f_l = -  E_{ 2l , 2 l+1 } +  E_{ 2 l+1 ,  l }  \]
 For $\alpha_0 $,
 \[ e_0 = E_{l+1, 1 } t , \; \; f_0 = E_{ 1 , l+1 } t^{-1} \]

Let
\begin{equation}
\Lambda^{\vee}  = \rm{Span}_{\Bbb Z} (  \alpha_0^\vee , \alpha_1^\vee , \dots , \alpha_l^\vee , 2d )
\end{equation}
 and $\Lambda $ be the dual lattice in $\eta^*$, the dual of $\eta$.
We have the root datum
\begin{equation}
   ( \Lambda^{\vee} , \Pi^\vee , \Lambda , \Pi )
\end{equation}
which is simply connected in the sense that $\Lambda^{\vee} / {\Bbb Z} \Pi^\vee $ is  torsion free as
${\Bbb Z}$-module.

\

The invariant symmetric bilinear form  (2.1) restricts to an invariant symmetric bilinear form on
 $ \hat{sl}_{2l+1}^\sigma $, it restricts further to a $W$-invariant symmetric bilinear form on $\eta$, where $W$ is still the Weyl group of the fixed point Lie algebra $\widehat{sl}_{2l+1}^\sigma$. 
We define $( \; , \; ):   \eta  \times  \eta \to {\Bbb Z}$ be $\frac 12 $ times this form. So we have
\[   ( \epsilon_i^* , \epsilon_j^* ) = \delta_{ij} ,   \; \; ( K , 2 d )  =1 ,  \; \;  ( K, K) = ( d , d ) = 0   .  \]
Its restriction on $ \Lambda^\vee$ is ${\Bbb Z}$-valued and defines a quadratic form on
 \[ Q:  \Lambda^\vee \to {\Bbb Z} ,  \; \;  Q ( x ) = ( x , x ) . \]

\

\noindent {\bf 3. Symplectic Realization of $A_{2l-1}^{(2)}$ ($l\geq 3$).}
In this case ${\frak g}  = sl_{2l} $.
The $\sigma $-fixed point algebra  $ {\frak g}_0 $ is the orthogonal Lie algebra of type $D_l$, it has
 Cartan subalgebra
\begin{equation}
 \eta_0 = \{  \rm{diag} ( d_1 , \dots , d_l , - d_1 , \dots , - d_l ) \}
 \end{equation}
 We write
\[ \rm{ diag} (d_1 , \dots , d_l , - d_1, \dots , - d_l ) = d_1 \epsilon_1 + \dots + d_l  \epsilon_l  ,\]
 so $\eta_0$ has a basis $ \epsilon_1 , \dots ,  \epsilon_l $. We define a symmetric bilinear form
  on $\eta_0$ by
 \[ ( \epsilon_i ,\epsilon_j ) = \delta_{ij} ,\]
   we identify the dual $\eta_0^*$ with $\eta $ by this form.
 The set of roots  of ${\frak g}_0$ is
\[ \Delta_0 \deff \{  \pm \epsilon_i \pm \epsilon_j ,    1\leq  i \ne j \leq l \}. \]
  The simple roots are
 \[ \Pi_0 = \{  \alpha_{1} =  \epsilon_1 - \epsilon_2 ,  \alpha_{2}=  \epsilon_2 - \epsilon_3, \dots , \alpha_{1-1} = \epsilon_{l-1} - \epsilon_l,  \alpha_l =  \epsilon_{l-1} + \epsilon_l \}\]
is a set of simple roots.  An element
$  \left( \begin{matrix}  a & b \\  c & d \end{matrix} \right)  \in sl_{2l} $ is in ${\frak g}_1$ iff
\[ a^T =  d  , \; \; \;   b - b^T = 0 , \; \; \;  c - c^T = 0   \]
The set of  weights of ${\frak g}_1$ is a highest weight ${\frak g}_0$-module with weights
\[\Delta_1  \deff  \{  \pm (   \epsilon_i  + \epsilon_j )  , 1 \leq  i , j \leq l ,\} \sqcup \{     \epsilon_i - \epsilon_j ,  1\leq i , j \leq l \}   \]
and $ 2 \epsilon_1 $ is the highest weight.

The twisted Lie algebra  $\widehat{sl}_{2l}^{\sigma } $ is of type $A_{2l-1}^{ ( 2 )} $.
It has  Cartan sublagebra
\begin{equation}
    \eta = \eta_0 +  {\Bbb C} K +  {\Bbb C} d
\end{equation}
It has a basis $\epsilon_1 , \dots, \epsilon_l , K , d $.
 Its real roots are
\[  \Phi_{re} = \{  \alpha + 2k \delta , \alpha \in   \Delta_0  ,k \in {\Bbb Z} \} \sqcup
    \{  \alpha + ( 2k -1 )\delta , \alpha \in   \Delta_1 - \{ 0 \} , k \in {\Bbb Z} \} .\]
 The simple roots  are
\begin{equation}
 \Pi = \{  \alpha_0 =  \delta - 2 \epsilon_1 , \;  \alpha_{1} , \dots , \alpha_l \}
\end{equation}
where $\alpha_1 , \dots , \alpha_l $ are simple roots for ${\frak g}_0$ as above.
The simple coroots are
\[  \Pi^\vee = \{ \alpha_0^{\vee} = K - \epsilon_1 ,  \alpha_{1}^{\vee} =  \epsilon_{1} - \epsilon_2 , \dots , \alpha_l^\vee =  \epsilon_{l-1} + \epsilon_l \} \]

\

We fix a $sl_2$-triple for each simple root $\alpha_i$.  For $ 1\leq  i \leq l-1$,
\[     e_i = E_{ i ,  i +1  } -  E_{ l + i + 1 , l + i    }  , \; \;   f_i =  E_{ i +1 , i  } -  E_{ l + i  ,  l + i +1    } . \]
For $\alpha_l = \epsilon_{l-1} + \epsilon_{l} $, we choose the
 \[  e_l = E_{ l-1 , 2 l } - E_{ l , 2 l-1 } , \; \; \;  f_l =  E_{ 2l , 1} - E_{ 2l-1 , 2 }  \]
For $  \alpha_i =  \epsilon_{ l-i } -\epsilon_{ l -i+1 } $, $ 1\leq i \leq l-1$,
    For $ \alpha_0 =   \delta - 2 \epsilon_1    $:
\[     e_0 =  E_{ l+1 , 1 } t  ,  \; \; \;   f_0 =  E_{ 1 , l+1 } t^{-1} \]
 One checks directly that
 \[   [ e_0 , f_0 ] =  \alpha_0^{\vee} =  - E_{1 , 1} +  E_{ l+1 , l+1 } + K = K - \epsilon_1 .\]
 The generators $ e_i , f_i , \alpha_i^\vee $,   $0 \leq  i \leq l $, are Chevalley generators.
  Let
  \[ \Lambda^\vee = \rm{Span}_{\Bbb Z} \{  \alpha_0^\vee , \dots , \alpha_l^\vee , 2 d \} \]
    $\Lambda $ be the dual lattice,
  We have
   root datum $ ( \Lambda^\vee , \Pi^\vee , \Lambda , \Pi )$. As in the case $A_{2l}^{(2)}$, we let
    $ Q: \Lambda^\vee \times \Lambda^\vee \to {\Bbb Z} $ be $\frac 12 $ times the restriction of (2.1) on
     $ \Lambda^\vee$.

\

\section{\Large Weil Representation.   }\label{Section 3}
\subsection{Weil Representations for Infinite Dimensional Symplectic Groups}
Let $F$ be a local field, we fix a non-trivial additive character $\psi : F \to U(1)$. It determines a unique  Haar measure $dy$  on $F$ such that
 the Fourier transform
  \[  f(x) \mapsto  \int_F f ( y ) \psi (x y ) d y \]
 defined using the measure $dy$ is an isometry on $ L^2 ( F , d y ) $.
  Using results in \cite{Z}, we define a projective representation of the twisted loop $ L G ( F ) $ defined in Section 1.3. Then we describe explicitly the action of $\rm{SL}_2 (F)$ corresponding to the
   simple root $\alpha_0$. Our formula shows that  the projective representation, when restricted to the above $\rm{SL}_2 ( F)$, gives a representation of its the two-fold metaplectic cover. 

We recall the basic settings in \cite{Z} Section 2. Let $V$ be a vector space over $F$ (not necessarily finite-dimensional), let $V^*$ be its linear dual space. We define a symplectic form on $X=V\oplus V^*$ given by $$\la v+v^*,w+w^* \ra=w^*(v)-v^*(w)$$ Let $ \rm{Sp} ( X, V^* )$ be the space of symplectic isomorphisms $g$ on $X$ such that $ V^*g $ and $V^*$ are commensurable.
We write $ g \in \rm{  Sp ( X )}$ as a matrix
\[  g = \left[ \begin{array}{cc}
 \alpha  &  \beta \\ \gamma & \delta \end{array} \right] \]
 where $\alpha : V\to V,  \beta : V\to V^*, \gamma : V^* \to V, \delta : V^* \to V^*$.
 The condition that $ V^*g $ and $V^*$ are commensurable is equivalent to that ${\rm dim} ( {\rm Im} \gamma ) < \infty $.

We need the following result from \cite{Z}, which is Proposition 2.3 in \cite{Z}. 

\begin{prop} If $ g = \left[ \begin{array}{cc}
 \alpha  &  \beta \\ \gamma & \delta \end{array} \right] \in \rm{Sp} ( X )$ is an element in
  $\rm{Sp}( X, V^* )$, we define
 an operator $T(g) $ on $ {\cal S} (V) $ by
\begin{equation}\label{action}
( T(g) f ) ( x ) = \int_{V_g}    S_{ g}( x + x^* )   f ( x \alpha + x^* \gamma ) d ( x^*\gamma  ), \end{equation}
where
\[ S_{g} ( x + x^*) = \psi \left( \frac 12 \la x\alpha , x \beta \ra +
          \frac 12 \la x^* \gamma , x^* \delta \ra + \la x^* \gamma , x \beta \ra \right) ;\]
The integration is for variable $ x^* \gamma $ over $V_g = {\rm Im } ( \gamma )$ with respect to a Haar measure.
The map  $ g \mapsto T(g)$ gives a projective representation of $\rm{Sp}( X, V^* )$ on $ {\cal S} ( V ) $.
\end{prop}
 In general there are no canonical choices of Haar measures on ${\rm Im}\gamma_g$. Nevertheless, for some special elements we choose measures that are most convenient.
  For example, for those $g\in\rm{Sp}(X,V^*)$ with ${\rm Im}\gamma_g=\{0\}$, then we can take the counting measure and define
  \begin{equation}
  (T(g)f)(x)=\psi\left(\frac{1}{2}\la x\alpha,x\beta \ra\right)f(x\alpha)
  \end{equation}
   Let
   \[ P=\{g\in\rm{Sp}(X,V^*):\mathrm{Im}\gamma_g=\{0\}\} \]
The following lemma, which is Lemma 2.5 in \cite{Z}, will be used in Section 5. 
\begin{lemma}[Lemma 2.5 in \cite{Z}]
For $p=\begin{pmatrix}
	\alpha_p & \beta_p \\ \gamma_p & \delta_p
\end{pmatrix}\in P$, $g=\begin{pmatrix}
	\alpha_g & \beta_g \\ \gamma_g & \delta_g
\end{pmatrix}\in \mathrm{Sp}(X,V^*)$,
 \begin{enumerate}
	\item We have $\mathrm{Im}(\gamma_g)=\mathrm{Im}(\gamma_{pg})$. If we choose the same Haar measures on $\mathrm{Im}(\gamma_g)$ and $\mathrm{Im}(\gamma_{pg})$, then $T(p)T(g)=T(pg)$, where $T(g)$ is the operator defined in (3.2.6) with respect to the Haar measure chosen above. 
	\item We have $\mathrm{Im}(\gamma_g)\alpha_p=\mathrm{Im}(\gamma_{gp})$. If we choose Haar measures on $\mathrm{Im}(\gamma_g)$ and $\mathrm{Im}(\gamma_{gp})$ so that $\alpha_p$ preserves the Haar measures, then $T(g)T(p)=T(gp)$.  
\end{enumerate}
In particular, for $p_1,p_2\in P$, we have $T(p_1)T(p_2)=T(p_1p_2)$. 
\end{lemma}

\par Now we apply the above general framework to the case of this paper. For $C$ as Section 2, we have a symplectic form on $ X = F(( t))^{n}$ defined by (1.1).
We have the Lagrangian space decomposition
\[     F(( t))^{n} =F [ t^{-1} ]^{n} t^{-1} \oplus F [[ t]]^{n}   = X_- \oplus X_+  .\] we take $V=X_-$, $V^*=X_+$. 
 Our twisted loop group $ L G ( F)  $ is a subgroup of $\rm{Sp}(X,X_+)$, by restriction, we have a projective representation of $LG(F)$ on ${\cal S} ( X_- )$,
 we call this representation the Weil representation of  $LG(F)$.

\

\subsection{The case of $A_{2l-1}^{(2)}$}
Let  $ X=F^{2l} (( t )) $ be the  symplectic space with symplectic form given in (1.1).
 Let $e_i \in F^{2n} $ denote the vector with $i$-th entry $1$ and all other entries $0$. Then
 $e_{i} t^j $ ($ 1 \leq i \leq 2l$, $ j \in {\Bbb Z} $) is a basis for $ F^{2l} [ t , t^{-1} ] $. We have
\[ \la e_{i} t^j ,   e_{ m } t^{ n } \ra = ( -1)^n \delta_{ i + m , 2 l+1 } \delta_{ j + n , -1 } .\]
\[  X =  X_+ \oplus X_- =  F^{2l} [[ t ]] \oplus F^{2l} [ t^{-1} ] t^{-1} .\]
	By Proposition 3.1 we have a projective representation of $ \rm{ Sp} ( X ,  X_+ ) $ on $ {\cal S} ( F^{2l} [ t^{-1} ] t^{-1} )$ given by $g\mapsto T(g)$. The twisted loop group
  given in (1.3)
   $$LG(F) =\{g(t)\in SL_{2l}(F((t))):C g ( -t )^{-T} C^{-1}  =  g ( t )\}$$
is a subgroup of $\Sp(X,X_+)$, hence acts projectively on ${\cal S} ( F^{2l} [ t^{-1} ] t^{-1} )$.

\par In the rest of this section, all the $2\times2$ matrices correspond to the $(1,l+1;1,l+1)$-submatrix of a $2l\times 2l$-matrix in $LG(F)$. We define the following elements in the $\rm{SL}_2$-subgroup corresponding to the root $\alpha_0$: \begin{equation}
	\begin{split}
		&X_0(r)=\begin{pmatrix}
	1 & \\ rt & 1
\end{pmatrix},\,X_{-0}(r)=\begin{pmatrix}
	1 & rt^{-1}\\  & 1
\end{pmatrix}\text{ for }r\in F,\\&W_0(r)=X_0(r)X_{-0}(-r^{-1})X_0(r),\,H_0(r)=W_0(r)W_0(1)^{-1}\text{ for }r\in F^*.
	\end{split}
\end{equation} Note that the root subgroups corresponding to other simple roots are contained in $P$, so the actions are simply defined by (3.2). We are going to define the operators $T(g)$ on $\cal S(X_-)$ for $g$ in (3.3).
 This is equivalent to choosing a Haar measure on $\rm{Im}\gamma_g$, but we use the similar approach as in \cite{Z}. Since our computation is only relavent to the $(1,l+1)$-variables of a Schwartz function, for $f\in\cal S(X_-)$ we write $f(\cdots,a_{-2},a_{-1}|a_0,a_1,\cdots)$ for the value of $f$ at $(\sum_{k<0}a_kt^k,\sum_{k\geq0}a_kt^{-1-k})$ (they stand for the $(1,l+1)$-variables. ) 
\subsubsection{Definition of $T(D_1)$ and $T(D_{-1})$} Before  defining operators $T(g)$ for $g$ in (3.3), we first define the operators $T(D_1)$ and $T(D_{-1})$ for the elements $$D_1=\rm{diag}(-t,1,\cdots,1,t^{-1},1,\cdots,1),\,\,D_{-1}=D_1^{-1}$$ where the $-t,t^{-1}$ lies in the $1,l+1$ diagonal entries respectively. Note that $D_1$ and $D_{-1}$ are not elements in the twisted loop group $LG(F)$, but they are elements in $\Sp(X,X_+)$, so the corresponding operators via Weil representation make sense. They will be used in the sequel when defining $T(X_{-0}(r))$. First we let $g=D_1$. For $(x,y)\in X_-=F^{2l}[t^{-1}]t^{-1}$ and $(u,v)\in X_+=F^{2l}[[t]]$ we have
 \begin{eqnarray*}
 &&   ( x  , y )   \alpha_g  = (  -( x t )_- ,  (y t^{-1})_- ) =  (  -( x t )_- ,  y t^{-1} )    \\
 &&   ( x  , y )   \beta_g  = (  -( x t )_+ ,  ( y t^{-1})_+ ) =  (  -( x t )_+ ,   0  )  \\
 &&   ( u   , v )   \gamma_g  = (  -( u t )_-   , ( v t^{-1} )_- )  = (   0  , ( v t^{-1} )_- ) \\
  &&   ( u   , v )   \delta_g  = (  -( u t)_+   , ( v t^{-1} )_+ )  = (   -u t   , ( v t^{-1} )_+ )
 \end{eqnarray*}
Here for a Laurent series $a(t)=\sum_n a_nt^n\in F((t))$, we define $a_-=\sum_{n<0}a_nt^n$ and $a_+=\sum_{n\geq0}a_nt^n$. So $ V_{g } =\{ ( 0 ,  a_0t^{-1} ):a_0\in F \}$ is obviously identified with $F$. We take the Haar measure on $V_g$ to be the Haar measure on $F$ under this identification. We see that
\[  \la  ( x  , y )   \alpha_g , ( x  , y )   \beta_g \ra = 0 , \; \; \;   \la ( u   , v )   \gamma_g ,  ( u   , v )   \delta_g \ra = 0 .\] Notice that for $ ( 0 , y_{-1} t^{-1}) \in V_g $,
 \[ \la   x^* \gamma ,   x \beta \ra =  \la ( 0 , y_{-1} t^{-1}  ) , - ( a_{-1} , 0 ) \ra=-a_{-1}y_{-1} \]
So by Proposition 3.1 we have
 \begin{eqnarray*}
T(D_1) f (    \dots ,   a_{-2} ,   a_{-1}  \, | \, a_0 , a_1 , a_2 , \dots )  =
 \int_{ F } f (    \dots ,  -a_{-3} , -a_{-2}  \, | \, y_{-1}, a_{0} , a_{ 1 } ,      \dots ) \psi ( -y_{-1} a_{-1}   )dy_{-1}
   \end{eqnarray*}
For $D_{-1}$ the situation is totally analogous to $D_1$ and we only record the result: \begin{eqnarray*}
    T(D_{-1}) f  (    \dots ,   a_{-2} ,   a_{-1}  \, | \, a_0 , a_1 , a_2 , \dots ) = \int_{ F } f (    \dots ,   -a_{-2} ,  -a_{-1} ,  y_0  \, | \,   a_{1} , a_{ 2 } ,      \dots )\cdot
  \psi (y_0 a_0) dy_0
   \end{eqnarray*}
   In particular, $T(D_1)$ and $T(D_{-1})$ are inverse of each other. 
\subsubsection{Definition of $T(X_0(r))$} Since $X_0(r)\in P$, the operator $T(X_0(r))$ is defined by (3.2), namely for $f\in\cal S(X_-)$ we have \begin{equation}
	(T(X_0(r))f)(\cdots,a_{-2},a_{-1}|a_0,a_1,\cdots)=\psi(\frac{1}{2}ra_0^2)f(\cdots,a_{-2}+ra_2,a_{-1}+ra_1|a_0,a_1,\cdots)
\end{equation}
\subsubsection{Definition of $T(X_{-0}(r))$} Next we consider $T(X_{-0}(r))$. Following \cite{Z} we define $T(X_{-0}(r))$ as follows: we have $X_{-0}(r)=D_{-1}\begin{pmatrix}
	1 & -rt\\&1
\end{pmatrix} D_1$ and $\begin{pmatrix}
	1 & -rt\\&1
\end{pmatrix}\in P$, so $T\begin{pmatrix}
	1 & -rt\\&1
\end{pmatrix}$ is defined by (3.2) and we simply define $$T(X_{-0}(r))=T(D_{-1})T\begin{pmatrix}
	1 & -rt\\&1
\end{pmatrix}T(D_1)$$
By (3.4) the operator $T\begin{pmatrix}
	1 & -rt\\ & 1
\end{pmatrix}$ is given explicitly as follows: \begin{align*}
	&(T\begin{pmatrix}
	1 & -rt\\ & 1
\end{pmatrix}f)(\cdots,a_{-2},a_{-1}|a_0,a_1,\cdots)\\=&\psi(-\frac12ra_{-1}^2)f(\cdots,a_{-2},a_{-1}|a_0-ra_{-2},a_1-ra_{-3},\cdots)
\end{align*}
\par Next we find the formula for $T(X_{-0}(r))$. For $f\in\cal S(X_-)$, we have \begin{equation}
	\begin{split}
		&(T(X_{-0}(r))f)(\cdots,a_{-2},a_{-1}|a_0,a_1,\cdots) \\=&(T(D_{-1})T\begin{pmatrix}
	1 & -rt\\&1
\end{pmatrix}T(D_1)f)(\cdots,a_{-2},a_{-1}|a_0,a_1,\cdots)\\=&\int_F(T\begin{pmatrix}
	1 & -rt\\&1
\end{pmatrix}T(D_1)f)(\cdots,-a_{-2},-a_{-1},-y|a_1,\cdots)\psi(a_0y)dy\\=&\int_F\psi(-\frac{1}{2}ry^2+a_0y)(T(D_1)f)(\cdots,-a_{-1},y|a_1+ra_{-1},a_2+ra_{-2},\cdots)dy\\=&\int_F\psi(-\frac{1}{2}ry^2+a_0y)(\mathcal{F}_0^{-1}f)(\cdots,a_{-2},a_{-1}|y,a_1+ra_{-1},\cdots)dy\\=&(\mathcal{F}_0\circ \psi(-\frac{1}{2}ry^2)\circ\mathcal{F}_0^{-1})f(\cdots,a_{-2},a_{-1}|a_0,a_1+ra_{-1},\cdots)
	\end{split}
\end{equation}
where $\psi(-\frac{1}{2}ry^2)$ is viewed as the multiplication operator on the $a_0$-variable and $\cal F_0$ is the Fourier transformation on the $a_0$-variable $$({\cal F}_0f) (  \dots ,  a_{-2} ,   a_{-1}  | a_0 , a_1 , a_2, \dots ) = \int_F    f (  \dots ,  a_{-1} |   y_0 , a_{1}, \dots  ) \psi (a_0 y_0 ) d y_0$$Since the conjugate of a multiplication operator by Fourier transform is the convolution operator by the Fourier transform of the function in the multiplication operator in the sense of distributions (see also \cite{W}), RHS of (3.5) is equal to $$(\mathcal{F}\psi(-\frac{1}{2}ry^2))*f(\cdots,a_{-2},a_{-1}|a_0,a_1+ra_{-1},\cdots)$$ $$=|r|^{-\frac{1}{2}}\gamma(-r)\int_F\psi(\frac{1}{2}r^{-1}(a_0-y)^2)f(\cdots,a_{-2},a_{-1}|y,a_1+ra_{-1},\cdots)dy$$ where $*$ is the convolution on the $a_0$-variable (in the sense of distributions), and $\gamma$ is the Weil index of $\psi$ determined by the following formula: $${\cal F}(\psi(\frac12cx^2))=\gamma(c)|c|^{-\frac12}\psi(-\frac12c^{-1}x^2)$$ Here we view the quadratic character $\psi(\frac12cx^2)$ as a tempered distribution on ${\cal S}(F)$ and $\cal F$ is the Fourier transform of tempered distributions. 
\par To summarize, we have \begin{equation}
	\begin{split}
		 &(T(X_{-0}(r))f)(\cdots,a_{-2},a_{-1}|a_0,a_1,\cdots)\\&=|r|^{-\frac{1}{2}}\gamma(-r)\int_F\psi(\frac{1}{2}r^{-1}(a_0-y)^2)f(\cdots,a_{-2},a_{-1}|y,a_1+ra_{-1},\cdots)dy
	\end{split}
\end{equation}
\subsubsection{Definition of $T(W_0(r))$} In the same spirit, we define $$T(W_0(r))=T(X_0(r))T(X_{-0}(-r^{-1})T(X_0(r))$$ and we can compute it by using formula (3.4) and (3.6) as follows: \begin{align*}
	 &(T(W_0(r))f)(\cdots,a_{-2},a_{-1}|a_0,a_1,\cdots)\\=&(T(X_0(r))T(X_{-0}(-r^{-1})T(X_0(r))f)(\cdots,a_{-2},a_{-1}|a_0,a_1,\cdots)\\=&\psi(\frac{1}{2}ra_0^2)(T(X_{-0}(-r^{-1}))T(X_0(r))f)(\cdots,a_{-2}+ra_2+,a_{-1}+ra_1|a_0,a_1,\cdots)\\=&|r|^\frac{1}{2}\gamma(r^{-1})\psi(\frac{1}{2}ra_0^2)\int_F\psi(\frac{1}{2}r(a_0-y)^2)(T(X_0(r))f)(\cdots,a_{-1}+ra_1|y,-r^{-1}a_{-1},-r^{-1}a_{-2},\cdots)dy\\=&|r|^\frac{1}{2}\gamma(r^{-1})\psi(\frac{1}{2}ra_0^2)\int_F\psi(-\frac{1}{2}r(a_0-y)^2+\frac{1}{2}ry^2))f(\cdots,ra_2,ra_1|y,-r^{-1}a_{-1},-r^{-1}a_{-2},\cdots)dy\\=&|r|^\frac{1}{2}\gamma(r^{-1})\int_F\psi(ra_0y)f(\cdots,ra_2,ra_1|y,-r^{-1}a_{-1},-r^{-1}a_{-2},\cdots)dy\\=&|r|^\frac{1}{2}\gamma(r^{-1})(\cal F_0f)(\cdots,ra_2,ra_1|ra_0,-r^{-1}a_{-1},-r^{-1}a_{-2},\cdots)dy
\end{align*} So we have \begin{equation}
	(T(W_0(r))f)(\cdots,a_{-2},a_{-1}|a_0,a_1,\cdots)=|r|^\frac{1}{2}\gamma(r^{-1})(T\begin{pmatrix}
		&-r^{-1}\\
		r&
	\end{pmatrix}\tau_1\mathcal{F}_0f)(\cdots,a_{-2},a_{-1}|a_0,a_1,\cdots)
\end{equation} where $\tau_1$ is the operator given by $$(\tau_1f)(\cdots,a_{-2},a_{-1}|a_0,a_1,\cdots)=f(\cdots,a_{-3},a_{-2}|a_{-1},a_0,\cdots)$$
\subsubsection{Definition of $T(H_0(r))$} Finally we define $T(H_0(r))$ by $$T(H_0(r))=T(W_0(r))T(W_0(1))^{-1}$$ Then by using (3.7) we can compute \begin{align*}
	&T(H_0(r))=T(W_0(r))T(W_0(1))^{-1}\\=&|r|^\frac{1}{2}\gamma(r^{-1})T\begin{pmatrix}
		&-r^{-1}\\
		r&
	\end{pmatrix}\tau_1\mathcal{F}_0\mathcal{F}_0^{-1}\tau_1^{-1}T\begin{pmatrix}
		&-1\\
		1&
	\end{pmatrix}^{-1}\gamma(1)^{-1}\\=&|r|^\frac{1}{2}\gamma(r^{-1})\gamma(1)^{-1}T\begin{pmatrix}
		r^{-1}&\\
		&r
	\end{pmatrix}
\end{align*}
So we have \begin{equation}
	T(H_0(r))=|r|^\frac{1}{2}\gamma(r^{-1})\gamma(1)^{-1}T\begin{pmatrix}
		r^{-1}&\\
		&r
	\end{pmatrix}
\end{equation}
In particular, we have \begin{align*}
	&T(H_0(r_1))T(H_0(r_2))T(H_0(r_1r_2))^{-1}\\=&\gamma(r_1^{-1})\gamma(r_2^{-1})\gamma (r_1^{-1}r_2^{-1})^{-1}\gamma(1)^{-1}=(r_1,r_2)_{Hil}
\end{align*} by the following property of Weil indices: $$\gamma(a)\gamma(b)=(a,b)_{Hil}\gamma(ab)\gamma(1)$$ Here $(r_1,r_2)_{Hil}$ is the quadratic Hilbert symbol over the field $F$. So the restriction of the projective representation $T$ to the $\mathrm{SL}_2$-subgroup corresponding to the simple root $\alpha_0$ can be lifted to a representation of the metaplectic $2$-fold cover of this $\mathrm{SL}_2(F)$. 

\subsection{The case of $A_{2l}^{(2)}$}
\par In this case the symplectic space is $X=F^{2l+1}((t))$ with symplectic form (1.1). We have the polarization
$$X=X_-\oplus X_+=F^{2l+1}[[t]]\oplus F^{2l+1}[t^{-1}]t^{-1}.$$
By Proposition 3.1,  $\mathcal{S}(F^{2l+1}[t^{-1}]t^{-1})$ admits a projective action of $\Sp(X,X_+)$.
The twisted loop group  $LG(F)$ is a subgroup of  $\Sp(X,X_+)$, it acts on  $\mathcal{S}(X_-)$ projectively by restriction.
\par In this section, all the $2\times2$ matrices correspond to the $(1,l+1;1,l+1)$-submatrix of a $(2l+1)\times(2l+1)$-matrix in $LG(F)$. We also define the elements in the $\rm{SL_2}$-subgroup corresponding to the root $\alpha_0$: \begin{equation}
	\begin{split}
		&X_0(r)=\begin{pmatrix}
	1 & \\ rt  & 1
\end{pmatrix},\,X_{-0}(r)=\begin{pmatrix}
	1 & rt^{-1}\\  & 1
\end{pmatrix}\text{ for }r\in F,\\&W_0(r)=X_0(r)X_{-0}(-r^{-1})X_0(r),\,H_0(r)=W_0(r)W_0(1)^{-1}\text{ for }r\in F^*.
	\end{split}
\end{equation} We define similarly the operators $T ( X_{- 0} ( r) ),  T ( X_{ 0} ( r) ), T ( W_{ 0} ( r) ), T ( H_0 ( r ) ) $ as in the previous case.
 The computations are almost identical to the previous case,  and the resulting formulas are identical to the formulas (3.4) (3.6) (3.7) (3.8).

\

\section{\Large Kac-Moody Group Functors and Metaplectic Covers }\label{Section 4}

In this section we first recall the Kac-Moody group functor defined by Tits \cite{Tits}
 and Patnaik-Puskas' construction of
  the metapectic covers of Kac-Moody groups \cite{PP}. We then give a description of the metaplectic cover in terms of generators and
   relations. In this section, all the rings are commutative with identity element $1$, and all functors are from the category of commutative unital rings to the category of groups.
\subsection{Root Data, Weyl Groups, and Lie Algebras.}
Recall that a \emph{root datum} is a quadruple
$D=\{\Lambda,\Lambda^\vee,\Pi=\{\alpha_i\}_{i\in I},\Pi^\vee=\{\alpha_i^\vee\}_{i\in I}\}$ where \begin{itemize}
	\item $\Lambda,\Lambda^\vee$ are free $\mathbb{Z}$-modules of the same rank, equipped with a perfect pairing $\la-,-\ra:\Lambda^\vee\times\Lambda\to\mathbb{Z}$ which identifies $\Lambda^\vee$ with the dual of $\Lambda$.
	\item $\{\alpha_i\}_{i\in I}\subseteq \Lambda$, $\{\alpha_i^\vee\}_{i\in I}\subseteq \Lambda^\vee$ are $\mathbb{Z}$-linearly independent.
	\item The matrix $A=(A_{ij}=\la \alpha_i^\vee,\alpha_j\ra)_{i,j\in I}$ is a \emph{generalized Cartan matrix}, namely $A_{ij}\in\mathbb{Z},A_{ii}=2,A_{ij}\leq0$ for $i\neq j$ and $A_{ij}=0\Leftrightarrow A_{ji}=0$.
\end{itemize}
A root datum will be denoted $(D,A)$ where $A$ is the generalized Cartan matrix. The sets $\Pi,\Pi^\vee$ are called \emph{simple roots, simple coroots} of the root datum $D$ respectively. From now on we fix a simply connected root datum $(D, A)$ in our discussion.
   The root datum $D$ is called simply connected if the ${\Bbb Z}$-module $\Lambda^\vee /  \rm{Span} (\Pi^\vee ) $ is torsion free. This definition
    is equivalent to Definition 2.1.6 in \cite{PP}.

 The Kac-Moody Lie Algebra $\mathfrak{g}(A)$ associated to $A$ is the Lie algebra over $\mathbb{C}$ generated by $e_i,f_i,\alpha_i^\vee$ for $i\in I$ subject to the following relations:
  \[ [\alpha_i^\vee, e_j]=A_{ij}e_j,\; [\alpha_i^\vee,f_j]=-A_{ij}f_j, \; [\alpha_i^\vee,\alpha_j^\vee]=0,\text{  }[e_i,f_i]=\alpha_i^\vee, \]
  \[ \text{for }i\neq j,\text{ }[e_i,f_j]=0,\text{  }(\ad e_i)^{-A_{ij}+1}(e_j)=(\ad f_i)^{-A_{ij}+1}=0. \]
 We denote the Weyl group by $W$ and the set of real roots by $\Phi_{re}$. The Weyl group $W$ has generators $s_i$, $ i \in I$.
   Let $\Phi_{re,+}$ ($\Phi_{re,-}$ be the set of positive (negative) real roots.
Let $s_i^*$ be the automorphism $\exp\ad e_i\cdot\exp\ad (-f_i) \cdot \exp\ad e_i$ of $\mathfrak{g}(A)$, let $W^*$ be the subgroup of $\mathrm{Aut}(\mathfrak{g}(A))$ generated by $s_i^*$ for $i\in I$, then we have a surjective homomorphism $v:W^*\to W$, $s_i^* \mapsto s_i$.
 For any real root $\alpha\in\Phi_{re}$, suppose $\alpha=w\alpha_i$ for $w\in W,i\in I$, let $w^*$ be a lifting of $w$ to $W^*$, then the set of opposite elements $w^*\{e_i,-e_i\}$ only depends on $\alpha$. We denote this set by $E_\alpha$, let $\mathfrak{g}_\alpha$ be the Lie subalgebra generated by $E_\alpha$.
\subsection{Tits Kac-Moody Group Functor}
Let $(D,A)$ be a root datum.
For every real root $\alpha$, we define $U_\alpha$ to be the unique group scheme over $\mathbb{Z}$ isomorphic to the additive group $\mathbb{G}_a$ whose Lie algebra is the $\mathbb{Z}$-subalgebra $\mathfrak{g}_{\alpha,\mathbb{Z}}$ generated by the set $E_\alpha$ over $\mathbb{Z}$. Every choice of $e_\alpha\in E_\alpha$ determines an isomorphism $x_\alpha:\mathbb{G}_a\to U_\alpha$. For a commutative ring $R$, we have $U_\alpha(R)=\{x_\alpha(r):r\in R\}$ and clearly $x_\alpha(r)x_\alpha(s)=x_\alpha(r+s)$ for $r,s\in R$.

A subset $\Psi\subseteq\Phi_{re}$ is called \emph{pre-nilpotent} if there exists $w,w'\in W$ such that $w\Psi\subseteq\Phi_{re,+}$ and $w'\Psi\subseteq\Phi_{re,-}$. It is called \emph{nilpotent} if it is complete, namely $\alpha, \beta\in \Psi,\alpha+\beta\in\Phi_{re}$ implies $\alpha+\beta\in\Psi$. Now let $\Psi\subseteq\Phi_{re}$ be a nilpotent set. The complex Lie algebra $\mathfrak{g}_\Psi=\oplus_{\alpha\in\Psi}\mathfrak{g}_\alpha$ is nilpotent, the lattice $\mathfrak{g}_\Psi=\oplus_{\alpha\in\Psi}  {\Bbb Z} E_\alpha$ is
 its ${\Bbb Z}$-form. It defines  a group scheme $U_\Psi$ over ${\Bbb Z}$
 together with embeddings $U_\alpha\hookrightarrow U_\Psi$ for $\alpha\in\Psi$ such that $U_\Psi(\mathbb{C})$ is the unique unipotent group corresponding to $\mathfrak{g}_\Psi$, and for any total ordering on $\Psi$ the resulting map $\prod_{\alpha\in\Psi}U_\alpha\to U_\Psi$ is an isomorphism of schemes.
In particular, for a pre-nilpotent pair of real roots $(\alpha,\beta)$, let $[\alpha,\beta]$ be the nilpotent set $\{m\alpha+n\beta:m,n\in\mathbb{Z}_{\geq0}\}\cap\Phi_{re}$, $]\alpha,\beta[=[\alpha,\beta]\backslash\{a,b\}$, then after choosing $e_\theta\in E_\theta$ for $\theta\in]\alpha,\beta[$ and a total ordering on $]\alpha,\beta[$ there exists well-defined integers $k(\alpha,\beta;\gamma)$ for $\gamma\in]\alpha,\beta[$ such that the following relation holds in $U_{[\alpha,\beta]}(R)$ for any commutative ring $R$:
 \begin{equation}\label{rel}
[x_\alpha(r),x_\beta(s)]=\prod_{\substack{\gamma=m\alpha+n\beta\\\gamma\in]\alpha,\beta[}}x_\gamma(k(\alpha,\beta;\gamma)r^ms^n),\text{ }\forall{r,s\in R}
\end{equation}

Now for every real root $\theta$, we fix a choice of $e_\theta\in E_\theta$, let $x_\theta:\mathbb{G}_a\to U_\theta$ be the associated isomorphism. The \emph{Steinberg group functor} $\St$ is defined by sending a ring $R$ the quotient of the free product of the groups $U_\alpha(R)$ for $\alpha\in\Phi_{re}$ by the smallest normal subgroup containing the relations
 (\ref{rel})
  for every pre-nilpotent pair $(\alpha,\beta)$ of real roots.

We define the following elements in $\St(R)$:
\[ w_\alpha(r)=x_\alpha(r)x_{-\alpha}(-r^{-1})x_\alpha(r),\text{ }h_\alpha(r)=w_\alpha(r)w_\alpha(-1)\text{ for }\alpha\in\Phi_{re},r\in R^* \]
and let
\[ \dot w_\alpha=w_\alpha(1),\, \dot w_i=\dot w_{\alpha_i},h_i(r)=h_{\alpha_i}(r).\]
We define the functor $\bold H$ by
\[ \bold H(R)=\mathrm{Hom}_\mathbb{Z}(\Lambda,R^*)\]
for a ring $R$. For $\lambda^\vee\in\Lambda^\vee$ and $r\in R^*$, write $r^{\lambda^\vee}\in\bold H(R)$ for the element sending $\mu\in\Lambda$ to $s^{\la\mu,\lambda^\vee\ra}$. The action of $W$ on $\Lambda$ induces a $W$-action on the functor $\bold H$. For example, for a simple reflection $s_i\in W$, we have $s_i(r^{\lambda^\vee})=r^{\lambda^\vee-\la\alpha_i,\lambda^\vee\ra\alpha_i^\vee}$.

The \emph{Tits Kac-Moody group functor} $\bold G_{D}$ (or usually $\bold G$ for short) is the functor sending a commutative ring $R$ to the quotient of the free product of $\St(R)$ and $\bold H(R)$ by the following relations: \begin{itemize}
	\item $tx_\alpha(r)t^{-1}=x_\alpha(t(\alpha)r)$ for $t\in\bold H(R),r\in R,\alpha\in\Phi_{re}$.
	\item $\dot w_i t\dot w_i^{-1}=s_i(t)$ for $i\in I,t\in\bold H(R)$.
	\item $h_i(r)=r^{\alpha_i^\vee}$ for $i\in I,r\in R^*$.
	\item $\dot w_ix_\alpha(r)\dot w_i^{-1}=x_{s_i\alpha}(\eta(\alpha_i,\beta)r)$ for $\alpha\in\Phi_{re},r\in R,i\in I$.
\end{itemize}
Here $\eta(\alpha_i.\beta)=\pm1$ is the sign coming from the ambiguity of choice of $e_\alpha\in E_\alpha$ for $\alpha\in\Phi_{re}$. It is determined by the equality $s_i^*(e_\alpha)=\eta(\alpha_i,\alpha)e_{s_i\alpha}$.

\subsection{Metaplectic Cover of Kac-Moody Group Functor}

Let $D$ be a root datum as before, $F$ be  a local field.  We denote the Kac-Moody group ${\bold G}_D ( F ) $ by $G$,
 let $U$ be the subgroup generated by $ x_a ( r )$, $a\in \Phi_{re, +} $, $r\in F$.
 For each $a \in \Pi $, we have a factorization $ U = U_a U^a $, where $U^a $ is generated by
the root subgroups $ U_b $ for $ b \in \Phi_{re, + } - \{ a \} $.

Let  $  ( , ) :  F^*\times F^* \to A $ be a Steinberg symbol with values in an abelian group $A$ in the sense of \cite{PP} (1.1.3).
 Let
 \begin{equation}
   Q : \Lambda^\vee \to  {\Bbb Z}
 \end{equation}
 be a $W$-invariant quadratic form and $B:  \Lambda^\vee \times \Lambda^\vee \to  {\Bbb Z} $ be
 the associated bilinear form given by
  \[ B ( x , y ) = Q ( x + y ) - Q( x ) - Q(y) .\]

 Since $D$ is simply connected,  we can extend $\Pi^\vee $ to a ${\Bbb Z}$-basis $\Pi^\vee_e$ of $\Lambda^\vee $.
  From the datum
   \[ D, \Pi^\vee_e , F, ( \; , \; ) ,  Q  \]
 Patnaik and Puskas constructed a metaplectic cover $ E$ of ${\rm \bf G}_D ( F)$ \cite{PP}.
  We recall the main steps in their construction. First a central extension
  $\widetilde H$ of $H$ is constructed using the Steinberg symbol $( \; )$ and $Q$.
  The group $\widetilde H$ is generated by $A$ and the symbols $\widetilde h_a(s)$ with $s\in F^*,\,a\in\Pi_e^\vee$
 subject to the following relations: $A$ is an abelian subgroup,  is in the center and
\begin{eqnarray}
&& 	\widetilde h_a(s)\widetilde h_a(t)\widetilde h_a(st)^{-1}=(s,t)^{Q(a^\vee)} \; \; {\rm for} \; a^\vee\in\Pi_e^\vee,  \; {\rm and} \; s,t\in F^*. \\
&& [\widetilde h_a(s),\widetilde h_b(t)]=(s,t)^{B(a^\vee,b^\vee)}\; {\rm for} \; a^\vee,b^\vee\in\Pi_e^\vee.
\end{eqnarray}
 We have the central extension
\[ 0\to A\to \widetilde H\xrightarrow{\varphi} H\to1 ,\]
 where the map $\varphi:\widetilde H\to H$ sends $\widetilde h_a(s)$ to $s^{a^\vee}\in H$ for each $a^\vee\in\Pi_e^\vee$ and $s\in F^*$.

\par Let $\widetilde H_{\mathbb Z}\subseteq\widetilde H$ be the subgroup generated by the elements $\widetilde h_a(-1)$ for $a\in\Pi$. For each $a\in\Pi$ we define a map $\frak s_a^{-1}:\widetilde H\to\widetilde H$ as $$\frak s_a^{-1}(x)=x\text{ for }x\in A,$$
$$\frak s_a^{-1}(\widetilde h_b(s))=\widetilde h_b(s)\widetilde h_a(s^{-\la a,b^\vee\ra})\text{ for }b^\vee\in\Pi_e^\vee,\,s\in F^*.$$
 The map $\frak s_a^{-1}$ is an automorphism of the group $\widetilde H$ whose inverse is denoted $\frak s_a$ (\cite{PP} Lemma 5.2.1 (2)).
 It has the properties that

\begin{prop}
	\begin{enumerate}
		\item  {\normalfont (\cite{PP} Lemma 5.2.1 (1))} For each $a\in\Pi$, we have $\frak s_a^{-1}(\widetilde h_a(s))=\widetilde h_a(s^{-1})$.
		\item {\normalfont (\cite{PP} Proposition 5.2.3)} The elements $\frak s_a^{-1} (a\in\Pi)$ satisfies the Braid relations.  Namely, $i\in I$ we set $\frak{s}_i^{-1}=\frak{s}_{\alpha_i}^{-1}$, then for $i\neq j\in I$ such that  $A_{ij}A_{ji}<4$ we have \begin{equation}
			\underbrace{\frak{s}_i^{-1}\frak{s}_j^{-1}\frak{s}_i^{-1}\cdots}_{h_{ij}}=\underbrace{\frak{s}_j^{-1}\frak{s}_i^{-1}\frak{s}_j^{-1}\cdots}_{h_{ij}}
		\end{equation} where $h_{ij}$ are integers defined by $$h_{ij}=\begin{cases}
			2 & A_{ij}A_{ji}=0;\\
			3 & A_{ij}A_{ji}=1;\\
			4 & A_{ij}A_{ji}=2;\\
			6 & A_{ij}A_{ji}=3;
		\end{cases}$$ (c.f. \cite{PP} (2.3)).
	\end{enumerate}
\end{prop}

Recall that $N$ is the subgroup of $G=\bold G_D(F)$ generated by $\{\dot w_a:a\in\Pi\}$ and $H=\bold H(F)$.
 In the next step, Patnaik-Puskas constructed a central extension $\tilde N $ of $N$
 $$1\to A\to\widetilde N\to N\to1.$$
The group $\tilde N $ is generated by $\tilde w_a $, $a\in \Pi $, $\tilde h_\gamma  ( s ) $ for $\gamma^\vee \in \Pi_e^\vee$.
 In the above map $\tilde N \to N $,
 $ \tilde w_a \mapsto \dot w_a$ and $\tilde h_\gamma  ( s ) \mapsto h_\gamma (s ) $.
 In $\tilde N$, we have relations
 \begin{equation}
    \tilde w_a \tilde h \tilde w_a^{-1} =  {\frak s}_a \tilde h
 \end{equation}
 for $\tilde h\in \tilde H $.

\

The metaplectic cover $E$ is constructed as a subgroup of permutations of certain fiber product set of $G$ with $\tilde N $.
It is generated by symbols $\lambda_a $ ($ a\in \Pi $), $\lambda (\tilde h ) $ ($\tilde h \in \tilde H$).

We have the following relations in $E$:
\begin{prop}{\normalfont (\cite{PP} Proposition 5.4.2)}
	\begin{enumerate}
		\item The map $\widetilde h\mapsto\lambda(\widetilde h)$ is an injective homomorphism $\lambda:\widetilde H\to E$ and $\lambda(A)$ is central in $E$.
		\item The map $u\mapsto\lambda(u)$ is an injective homomorphism $\lambda:U\to E$.
		\item For $u\in U$ and $\widetilde h\in\widetilde H$ we have $\lambda(\widetilde h) \lambda(u)\lambda(\widetilde h^{-1})=\lambda(\varphi(\widetilde h)u\varphi(\widetilde h)^{-1})$.
		\item For $a\in \Pi $, $u\in U^a$,  we have $\lambda_a\lambda(u)\lambda_a^{-1}=\lambda(\dot w_au\dot w_a^{-1})$ and also $\lambda_a^{-1}\lambda(u)\lambda_a=\lambda(\dot w_a^{-1}u\dot w_a)$.
		\item $\lambda_a^2=\lambda(\widetilde h_a(-1))$.
		\item $\lambda_a^{-1}\lambda(\widetilde h)\lambda_a=\lambda({\frak s}_a^{-1} \widetilde h )$.
		\item The elements $\lambda_a$ satisfy the braid relations (4.5).
	\end{enumerate}
\end{prop}
There is an exact sequence of groups
\[   1 \to  A\to E \to G \to 1 \]
where the surjective homomorphism $p:  E \to G$ have the properties that
$ p ( \lambda ( u ) ) = u $  for $ u \in U $,  $ p ( \lambda ( \tilde h  ) )= \phi ( h ) $,
$p ( \lambda_a ) = \dot w_a $.

\

To show that our projective presentation of loop group $LG(F) $ in Section 3 leads to a representation of
 the metaplectic group $E$ corresponding to $LG(F) $, we need to represent $E$ in terms of generators and relations.

  First we derive more relations among generators of $E$.
 We write $x_a ( r ) $ for  $\lambda ( x_a ( r ) )$ for $ a\in \Phi_{ re, +} $.
For each $ a\in \Pi$, we introduce the
\begin{equation}
   x_{-a} ( -r ) =  \lambda_a \lambda ( x_{a} ( r ) ) \lambda_a^{-1}
\end{equation}
This definition is motivated by (5.134) in \cite{PP}.
Then we introduce for $a\in \Pi $, $ r \in F^* $,
\begin{equation}
  w_a ( r ) = x_{a} ( r ) x_{-a} ( - r^{-1} ) x_a (  r )
\end{equation}

\begin{prop}
	In the group $E$ we have
$$w_a(-1)=\lambda_a$$
 where $w_a(r)$ is defined via formula (8) above.
\end{prop}
\paragraph{Proof.} We will used notations in \cite{PP} freely in the proof.
  By Lemma 5.4.1 in \cite{PP}, the group $E$ acts simply transitively on the set $S$, so it suffices to compute the action of $w_a(-1)$ and $\lambda_a$ on $\bold1_S=(1,1)\in S$ and compare the results. From the formula (6) of the action of $\lambda_a$ on $S$, we immediately have the following formula: $$\lambda_a^{-1}(g,\widetilde n)=\begin{cases}
		(\dot w_a^{-1}g,\widetilde w_a^{-1}\widetilde n), &\text{if } \nu(\dot w_a^{-1}g)=\dot w_a^{-1}\nu(g)\\
		(\dot w_a^{-1}g,\widetilde h_a(s)^{-1}\widetilde n), &\text{if } \nu(\dot w_a^{-1}g)=h_a(s)^{-1}\nu(g)
	\end{cases}$$
	We can compute as follows: \begin{align*}
		 &\lambda(w_a(r))\cdot(1,1)=\lambda(x_a(r))\lambda_a\lambda(x_a(r^{-1}))\lambda_a^{-1}\lambda(x_a(r))\cdot(1,1)\\&=\lambda(x_a(r))\lambda_a\lambda(x_a(r^{-1}))\lambda_a^{-1}\cdot(x_a(r),1)
	\end{align*} Since $\nu(\dot w_a^{-1}x_a(r))=\dot w_a^{-1}=\dot w_a^{-1}\nu(x_a(r))$, we have \begin{align*}
		 &\lambda(x_a(r))\lambda_a\lambda(x_a(r^{-1}))\lambda_a^{-1}\cdot(x_a(r),1)\\&=\lambda(x_a(r))\lambda_a\lambda(x_a(r^{-1}))\cdot(\dot w_a^{-1}x_a(r),\widetilde w_a^{-1})\\&=\lambda(x_a(r))\lambda_a\cdot(x_a(r^{-1})\dot w_a^{-1}x_a(r),\widetilde w_a^{-1})
	\end{align*} Since $\nu(\dot w_ax_a(r^{-1})\dot w_a^{-1}x_a(r))=\nu(x_{-a}(-r^{-1}))=h_a(-r)\dot w_a=h_a(r)\dot w_a^{-1}=h_a(r)\nu(x_a(r^{-1})\dot w_a^{-1}x_a(r))$, we have \begin{align*}
		&\lambda(x_a(r))\lambda_a\cdot(x_a(r^{-1})\dot w_a^{-1}x_a(r),\widetilde w_a^{-1})=\lambda(x_a(r))\cdot(x_{-a}(-r^{-1})x_a(r),\widetilde h_a(r)\widetilde w_a^{-1})\\&=(w_a(r),\widetilde h_a(r)\widetilde w_a^{-1})
	\end{align*} So we have $\lambda(w_a(-1))\cdot(1,1)=(w_a(-1),\widetilde h_a(-1)\widetilde w_a^{-1})$. On the other hand we obviously have $\lambda_a(1,1)=(\dot w_a,\widetilde w_a)$, so it suffices to prove that $\widetilde h_a(-1)\widetilde w_a^{-1}=\widetilde w_a$. This follows from the definition of $\widetilde N$: $\widetilde N$ is the quotient of $\widetilde N_\mathbb{Z}\ltimes \widetilde H$ by an ideal $J$ generated by elements of the form $(\widetilde w_a^{-2},\widetilde h_a(-1))$, which means that in $\widetilde N$ we have $\widetilde h_a(-1)=\widetilde w_a^2$, hence the result.

\begin{prop}
	In the  group $E$ we have $$w_a(r)w_a(-1)=\widetilde h_a(r)$$
\end{prop}
\paragraph{Proof.} We also compute the action on $\bold1_S=(1,1)\in S$. From the proof of the last proposition we have $\lambda(w_a(-1))\cdot(1,1)=(\dot w_a,\widetilde w_a)$ So we have \begin{align*}
	&\lambda(w_a(r)w_a(-1))\cdot(1,1)=\lambda(x_a(r))\lambda_a\lambda(x_a(r^{-1}))\lambda_a^{-1}\lambda(x_a(r))\cdot(\dot w_a,\widetilde w_a)\\&=\lambda(x_a(r))\lambda_a\lambda(x_a(r^{-1}))\lambda_a^{-1}\cdot(x_a(r)\dot w_a,\widetilde w_a)
\end{align*} Since $\nu(\dot w_a^{-1}x_a(r)\dot w_a)=\nu(x_{-a}(-r))=h_a(-r)^{-1}\dot w_a=h_a(-r)^{-1}\nu(x_a(r)\dot w_a)$, we have \begin{align*}
	&\lambda(x_a(r))\lambda_a\lambda(x_a(r^{-1}))\lambda_a^{-1}\cdot(x_a(r)\dot w_a,\widetilde w_a)\\&=\lambda(x_a(r))\lambda_a\lambda(x_a(r^{-1}))\cdot(x_{-a}(-r),\widetilde h_a(-r)^{-1}\widetilde w_a)\\&=\lambda(x_a(r))\lambda_a\cdot(x_a(r^{-1})x_{-a}(-r),\widetilde h_a(-r)^{-1}\widetilde w_a)
\end{align*} Since $\nu(\dot w_ax_a(r^{-1})x_{-a}(-r))=\nu(x_{-a}(-r^{-1})x_a(r)\dot w_a)=\nu(x_a(-r)h_a(r))=h_a(r)$ and $\nu(x_a(r^{-1})x_{-a}(-r))=h_a(-r^{-1})\dot w_a=\dot w_ah_a(-r)$, we have $$\nu(\dot w_ax_a(r^{-1})x_{-a}(-r))=\dot w_a\nu(x_a(r^{-1})x_{-a}(-r))$$ Thus \begin{align*}
	&\lambda(x_a(r))\lambda_a\cdot(x_a(r^{-1})x_{-a}(-r),\widetilde h_a(-r)^{-1}\widetilde w_a)\\&=\lambda(x_a(r))\cdot(\dot w_ax_a(r^{-1})x_{-a}(-r),\widetilde w_a\widetilde h_a(-r)^{-1}\widetilde w_a)\\&=(h_a(r),\widetilde w_a\widetilde h_a(-r)^{-1}\widetilde w_a)
\end{align*}By Proposition 4.1 (i) and (4.6) we have $\widetilde w_a\widetilde h_a(-r)^{-1}\widetilde w_a^{-1}=\widetilde h_a(-r^{-1})^{-1}$, so we have $$\lambda(w_a(r)w_a(-1))\cdot(1,1)=(h_a(r),\widetilde h_a(-1)\widetilde h_a(-r^{-1})^{-1})$$ On the other hand $$\lambda(\widetilde h_a(r))(1,1)=(h_a(r),\widetilde h_a(r))$$ By the properties in \cite{PP} 1.1.3 of a bilinear Steinberg symbol, we have $(r,-r^{-1})=(r,-1)(r,r^{-1})=(r,-1)^2=1$, thus $$\widetilde h_a(r)\widetilde h_a(-r^{-1})=\widetilde h_a(-1)(r,-r^{-1})^{\mathsf Q(a^\vee)}=\widetilde h_a(-1)$$ That is to say $\widetilde h_a(-1)\widetilde h_a(-r^{-1})^{-1}=\widetilde h_a(r)$. Hence the result.

\subsection{Generators and Relations for Mateplectic Cover $E$}

 \begin{theorem}
   The group $E$ is isomorphic to the group generated by symbols
  $\{\lambda (u):u\in U\}$, $\{\lambda(\widetilde h):\widetilde h\in \widetilde H\}$ and $\{\lambda_a\,(a\in \Pi)\}$
 with following relations
\begin{enumerate}
	\item For $A\subseteq \widetilde H$, $\lambda (A)\subseteq E$ is central in $E$.
	\item For $u\in U$ and $\widetilde h\in\widetilde H$,
    $\lambda (\widetilde h) \lambda (u)\lambda (\widetilde h^{-1})=\lambda (\varphi(\widetilde h)u\varphi(\widetilde h)^{-1})$.
	\item For $u\in U^a$,  $\lambda_a \lambda (u)\lambda_a^{-1}=\lambda (\dot w_au\dot w_a^{-1})$ and also $\lambda_a^{-1}\lambda(u)\lambda_a =\lambda (\dot w_a^{-1}u\dot w_a)$.
	\item $\lambda_a ^2=\lambda (\widetilde h_a(-1))$.
	\item $\lambda_a ^{-1}\lambda (\widetilde h)\lambda_a =\lambda({\frak s}_a^{-1} \widetilde h )$.
	\item The elements $\lambda_a $ satisfy the braid relations (4.5).
	\item Let $w_a (r):=\lambda (x_a(r))\lambda_a \lambda (x_a(r^{-1}))\lambda_a ^{-1}\lambda (x_a(r))$, then $w_a (-1)=\lambda_a$.
	\item $w_a (r)w_a(-1)=\lambda (\widetilde h_a(r))$.
\end{enumerate}
 \end{theorem}

 We denote $E'$ the group generated by $\{\lambda (u):u\in U\}$, $\{\lambda  (\widetilde h):\widetilde h\in \widetilde H\}$ and $\{\lambda_a\,(a\in \Pi)\}$ with the relations (1)-(8). To distinguish the generators of $E'$ and $E$, we write the generators
 $E'$  by $\{\lambda' (u):u\in U\}$, $\{\lambda' (\widetilde h):\widetilde h\in \widetilde H\}$ and $\{\lambda'_a\,(a\in \Pi)\}$,

  \par We define a homomorphism $\theta:E'\to E$ by sending the generators $\lambda'(\widetilde h),\,\lambda'(u),\,\lambda'_a$ in $E'$ to $\lambda(\widetilde h),\lambda(u),\lambda_a$ in $E$ respectively. This extends to a well-defined homomorphism because the relations (1)-(8) above all holds after mapping to $E$ via $\theta$: (1)-(6) are sent to the relations (1)-(7) in Proposition 4.2 , (7) and (8) are sent to the relations verified in the previous section. Clearly $\theta$ is a surjection. In particular, the homomorphism $\theta\circ\lambda':U\to E,\,u\mapsto\theta(\lambda'(u))$ is just the map $u\mapsto \lambda(u)$, which is injective by \cite{PP} prop 5.4.2 (2). So $\lambda':U\to E'$ is injective. Similarly $\lambda':\widetilde H\to E'$ is also injective.
\par To show that $\theta$ is injective we need the following combinatorical lemma, which is Corollary 2 in \cite{Bo} Chapter 6 Section 1.6: 
\begin{lemma}
	For $w\in W$, $a_i\in\Pi$, if $wa_i>0$ is a positive root, then for any reduced expression $w=s_{i_1}\cdots s_{i_k}$ we have $$s_{i_k}a_i,s_{i_{k-1}}s_{i_k}a_i,\cdots,s_{i_1}s_{i_2}\cdots s_{i_k}a_i$$ are all positive roots.
\end{lemma}
\paragraph{Proof of this combinatorical lemma:}Since $w(a_i)>0$ we have $\ell(ws_i)=\ell(w)+1$. If there exists some $1\leq p\leq k$ s.t. $s_{i_{p+1}}\cdots s_{i_k}a_i<0$, then $\ell(s_{i_{p+1}}\cdots s_{i_k}s_i)=\ell(s_{i_{p+1}}\cdots s_{i_k})-1=k-p-1$, thus $k+1=\ell(w)+1=\ell(ws_i)\leq\ell(s_{i_1}\cdots s_{i_{p}})+\ell(s_{i_{p+1}}\cdots s_{i_k}s_i)=p+k-p-1=k-1$, contradiction.

Since $ \lambda'_a $'s satisfy relations (4) (6) in Proposition 4.5, for every $ w \in W$, we can define $\lambda'_w$, which is unique to an element in $\lambda'(\widetilde H)$.

\begin{theorem}$E'$ admits the following decomposition:
	$$E'=\bigcup_{w\in W}\lambda'(U)\lambda'(\widetilde H)\lambda_w'\lambda'(U)$$
\end{theorem}
\paragraph{Proof.} The RHS set contains all the generators of $E'$, so it suffices to prove that RHS is closed under multiplication. For a generic element $$g=\lambda'(u_1)\lambda'(\widetilde h)\lambda'_w\lambda'(u_2)\text{ with }u_1,u_2\in U,\widetilde h\in\widetilde H,w\in W$$ in the RHS, it suffices to prove that left and right multiplication of this element by all three types of generator $\lambda'(u_0)\,(u_0\in U),\lambda'(\widetilde h_0)\,(\widetilde h_0\in\widetilde H),\lambda'_a\,(a\in \Pi)$ is still an element in the RHS. To do this we mimic the treatment of a BN-pair, let $\lambda'(\widetilde B)$ be  the subgroup of $E'$ generated by $\lambda'(U)$ and $\lambda'(\widetilde H)$ (by relation (E'2), $\lambda'(\widetilde H)$ normalizes $\lambda'(U)$, $\lambda'(\widetilde B)$ is a semi-direct product of $\lambda'(U)$ and $\lambda'(\widetilde H)$). Then we have $$\lambda'(U)\lambda'(\widetilde H)\lambda_w'\lambda'(U)=\lambda'(\widetilde B)\lambda'_w\lambda'(\widetilde B)$$ since $\lambda'_w$ normalizes $\lambda'(\widetilde H)$ by relation (E'5). Not it is clear that the RHS set $$\bigcup_{w\in W} \lambda'(\widetilde B)\lambda'_w\lambda'(\widetilde B)$$ is closed under left and right multiplication of elements in $\lambda'(\widetilde B)$, which contains all elements in $\lambda'(U)$ and $\lambda'(\widetilde H)$ (indeed, each "Bruhat cell" $\lambda'(\widetilde B)\lambda'_w\lambda'(\widetilde B)$ is closed under left and right multiplication by elements in $\lambda'(\widetilde B)$).
\par It suffices to prove that RHS is closed under multiplication by $\lambda_a$ for $a\in\Pi$, which follows from the assertion $$\lambda'_a\lambda'(\widetilde B)\lambda'_w\subseteq(\lambda'(\widetilde B)\lambda'_{w_aw}\lambda'(\widetilde B))\bigcup(\lambda'(\widetilde B)\lambda'_{w}\lambda'(\widetilde B))$$ Recall that we have a decomposition $U=U^aU_a$ and $\lambda'_a$ normalizes $\lambda'(U^a)$ by relation (E'3), so it suffices to prove $$\lambda'_a\lambda'(x_a(r))\lambda'_w\in(\lambda'(\widetilde B)\lambda'_{w_aw}\lambda'(\widetilde B))\bigcup(\lambda'(\widetilde B)\lambda'_{w}\lambda'(\widetilde B))$$ for any $r\in F$.
\par If $w^{-1}a>0$, by lemma 4.1 and relation (E'3) we have $\lambda'^{-1}_w\lambda'(x_a(r))\lambda'_w\in\lambda'(U)$. Also $w^{-1}a>0$ implies $\ell(w_aw)=\ell(w)+1$, thus $\lambda'_a\lambda'_w=\lambda'_{w_aw}$. So we have $$\lambda'_a\lambda'(x_a(r))\lambda'_w= \lambda'_a\lambda'_w(\lambda'^{-1}_w\lambda'(x_a(r))\lambda'_w)\in\lambda_{w_aw}\lambda'(U)\subseteq \lambda'(\widetilde B)\lambda'_{w_aw}\lambda'(\widetilde B) $$
\par The rest is to deal with the more sophisticated case $w^{-1}a<0$. To do this we first need
\begin{lemma}
	$$\lambda'(\widetilde B)\bigcup(\lambda'(\widetilde B)\lambda'_a\lambda'(\widetilde B))$$ is a subgroup of $E'$.
\end{lemma}
\paragraph{Proof of the lemma:}After using the decomposition $U=U^aU_a$ as above, it suffices to prove $$\lambda'_a\lambda'(x_a(s))\lambda'_a\in \lambda'(\widetilde B)\bigcup(\lambda'(\widetilde B)\lambda'_a\lambda'(\widetilde B))$$ By relation (E'4) $\lambda'^2_a=\lambda'(\widetilde h_a(-1))\in\lambda'(\widetilde H)$, so it suffices to prove $$\lambda'_a\lambda'(x_a(s))\lambda'^{-1}_a\in \lambda'(\widetilde B)\bigcup(\lambda'(\widetilde B)\lambda'_a\lambda'(\widetilde B))$$ If $s=0$ this is obvious, so in the followin we assume $s\in F^*$. By relation (E'8) and the definition of $w'_a(r)$ we have $$\lambda'(\widetilde h_a(r))\lambda_a^{-1}=w_a'(r)=\lambda'(x_a(r))\lambda_a'\lambda'(x_a(r^{-1}))\lambda_a'^{-1}\lambda'(x_a(r))$$ Thus $$\lambda'_a\lambda'(x_a(s))\lambda'^{-1}_a=\lambda'(x_a(-s^{-1}))\lambda'(\widetilde h_a(s^{-1}))\lambda'(\widetilde h_a(-1))\lambda'_a \lambda'(x_a(-s^{-1}))\in \lambda'(\widetilde B)\lambda'_a\lambda'(\widetilde B) $$ (note that this is an analogue of lemma 2.9 in the group $E'$) The lemma is proved.
\par We come back to the proof of the theorem. By the lemma above, $$\lambda'(\widetilde B)\bigcup(\lambda'(\widetilde B)\lambda'_a\lambda'(\widetilde B))$$ is a subgroup of $E'$, so we have $$\lambda'_a\lambda'(\widetilde B)\subseteq(\lambda'(\widetilde B)\lambda'_a)\bigcup(\lambda'(\widetilde B)\lambda'_a\lambda'(\widetilde B)\lambda'_a)$$ thus $$\lambda'_a\lambda'(\widetilde B)\lambda'_w\subseteq(\lambda'(\widetilde B)\lambda'_a\lambda'_w)\bigcup(\lambda'(\widetilde B)\lambda'_a\lambda'(\widetilde B)\lambda'_a\lambda'_w)$$ Since $w^{-1}a<0$, we have $\ell(w_aw)=\ell(w)-1$, so $\lambda'_a\lambda'_{w_aw}=\lambda'_w$, namely $\lambda'_a\lambda'_w=\lambda'(\widetilde h_a(-1))\lambda_{w_aw}$. So we can modify RHS of the above assertion, which gives $$\lambda'_a\lambda'(\widetilde B)\lambda'_w\subseteq(\lambda'(\widetilde B)\lambda'_{w_aw})\bigcup(\lambda'(\widetilde B)\lambda'_a\lambda'(\widetilde B)\lambda'_{w_aw})$$ Note that $(w_aw)^{-1}a>0$, so by the first case we already obtained, we have $$\lambda'_a\lambda'(\widetilde B)\lambda'_{w_aw}\subseteq(\lambda'(\widetilde B)\lambda'_{w}\lambda'(\widetilde B))\bigcup(\lambda'(\widetilde B)\lambda'_{w_aw}\lambda'(\widetilde B))$$ Plugin this to the second bracket of the above yields $$\lambda'_a\lambda'(\widetilde B)\lambda'_w\subseteq(\lambda'(\widetilde B)\lambda'_{w_aw}\lambda'(\widetilde B))\bigcup(\lambda'(\widetilde B)\lambda'_w\lambda'(\widetilde B))$$ which finishes the proof of the theorem.
\begin{theorem}
	$\theta:E'\to E$ is an isomorphism.
\end{theorem}
\paragraph{Proof.} Clearly $\theta$ is surjective. For injectivity of $\theta$, we take $k\in\ker\theta$ and try to prove $k=1$ in $E'$. By the above decomposition, we may assume $k=\lambda'(u_1)\lambda'(\widetilde h)\lambda'_w\lambda'(u_2)$ for some $u_1,u_2\in U,\,w\in W,\,\widetilde h\in\widetilde H$. Then $\theta(k)=\lambda(u_1)\lambda(\widetilde h)\lambda_w\lambda(u_2)=1$ in $E$. We look at the action of this element on $(1,1)\in S$. After this action the first factor becomes $u_1\varphi(\widetilde h)\dot wu_2\in G$, so $u_1\varphi(\widetilde h)\dot wu_2=1$ in $G$. By the Bruhat decomposition in $G$, this means $w=1,\,\varphi(h)=1,\,u_1u_2=1$. So $\lambda'_w=1$, $\widetilde h\in A$, which means $\lambda'(\widetilde h)$ is central, so $k=\lambda'(\widetilde h)\lambda'(u_1)\lambda'(u_2)=\lambda'(\widetilde h)$, and $\theta(k)=\lambda(\widetilde h)=1$. Since $\lambda:\widetilde H\to E$ is injective, this means that $\widetilde h=1$ in $\widetilde H$. Thus $k=1$ in $E'$.

\

\section{\Large Weil Representation of  Metaplectic Kac-Moody Group of Type $A_n^{(2)}$}\label{Section 5}

Let $F$ be a local field of characteristic zero. For the root datum $D$ as in Section 2 of type $A^{(2)}_n $. Apply the Tits functor to $D$,
  we have affine Kac-Moody group $G_D(F)$. And also we have the twisted loop group $LG (F) = SL_{n+1}(F[t , t^{-1}])^\sigma $ as in Section 1. We have surjective homomorphism
\begin{equation}  p:   G_D(F) \to  LG (F)  .\end{equation}
  which is defined as follows.
 The Kac-Moody Lie algebra ${\frak g} (A) $ defined in Section 2 can be defined on $F$, we denote this Lie algebra by
 ${\frak g} (A)_F$.
   As in Section 2, it has a realization as a central extension of $ sl_{n+1} ( F [ t , t^{-1} ] )^\sigma $.
  The space $ F [ t , t^{-1} ] ^{n+1} $ has an obvious $ sl_{n+1} ( F [ t , t^{-1} ] )^\sigma $-module structure.
   So it is a module over ${\frak g} (A)_F$ via the Lie algebra morphism
 \[ p :   {\frak g} (A)_F \to sl_{n+1} ( F [ t , t^{-1} ] )^\sigma .\]
 which induce the group homomorphism $p$. We describe the images of some elements.
For each real root,  $ e_\alpha $ acts as a nilpotent operator, so
  $ x_\alpha ( r ) = \sum_{ n=0}^\infty \frac 1 { n !}  e_\alpha^n$ is a well-defined operator on $ F [ t , t^{-1} ] ^{n+1} $,
   and this operator preserves the symplectic form and is in $LG(F)$, this gives the image of $ x_\alpha (r )$ in $ LG ( F) $.
   For $r^{\lambda^\vee} \in  {\bf H} ( F) $,  if $\lambda^\vee \in \rm{Span} ( \alpha_i^\vee , i \geq 0 ) $,
   it acts on $ F [ t , t^{-1} ] ^{n+1} $ in an obvious way, so it gives element in $LG( F)$.
 For $\lambda^\vee = 2 d $,
 \[  p ( r^{2d} ) v ( t) = v ( r^2 t ).  \]

We take the Steinberg symbol $(  \, , \, ) :  F^*\times F^* \to \{ \pm 1 \}$ in Section 4 as the Hilbert symbol,
 we have the datum $ ( D,  \Pi^\vee_e,  F , ( \, , \, ) , Q ) $, where $\Pi^\vee_e $ and $Q$ are as in Section 2.
Let $E$ be the metaplectic group constructed by Patnaik-Puskas as in Section 4. We denote the image of $g \in E $ in $LG( F) $ under the composition
  $E \to G_D ( F) $ ,  $ p : G_D (F) \to LG ( F) $ by $ \bar g $.

 \begin{theorem}
   The projective Weil representation of $LG(F)$ on ${\cal S} ( F [ t^{-1} ]^{n+1} t^{-1} ) $
   lifts to a representation of $E$ on ${\cal S} ( F [ t^{-1} ]^{n+1} t^{-1} ) $.
  \end{theorem}

\noindent {\it Proof.}  By Theorem 4.5, we need to define operators
   $\pi ( \lambda (u))$, $ u\in U$, $\pi( \lambda(\widetilde h) ) $, $\widetilde h\in \widetilde H $ and $\pi ( \lambda_a)$, $a\in \Pi$
on ${\cal S} ( F [ t^{-1} ]^{n+1} t^{-1} ) $ that satisfy the relations (1)-(8) in Theorem 4.5.
 For $ \lambda (u)(u\in U)$,  we have $ \overline {\lambda ( u) } \in P$ (recall that $P$ is defined before Lemma 3.2), and
  we define $ \pi ( \lambda (u))= T  ( \overline {\lambda ( u) } ) $.
  For $ \alpha_i \in \Pi^\vee $, $ 1\leq  i \leq l $, $s\in F^*$, we also have $\overline{h_{\alpha_i}(s)}\in P$, so we define $\pi(\tilde h_{\alpha_i} ( s ) )  = T( \overline{  \tilde h_{\alpha_i} (s ) } ) $. For $a\in A\subseteq E$ we define $\pi(\lambda(a))=a$ (note that in this case $A=\{\pm1\}$). We define $ \pi ( \tilde h_{\alpha_0} ( s ) )$ to be the operator $T(H_0(s))$ defined in (3.8). $\pi ( \tilde h_{2d} ( s ) )$ acts on $ {\cal S} ( F [ t^{-1} ]^{n+1} t^{-1} ) $ by
  \[   \pi ( \tilde h_{2d} ( s ) ) f ( v ( t) ) = f ( v ( s^2 t ) ). \]
 For $\pi ( \lambda_{\alpha_i } )$, $ 1\leq i \leq l $,
  $ \overline{ \lambda_{\alpha_i} } \in P $, so we define
  $\pi ( \lambda_{\alpha_i } ) = \pi ( \overline { \lambda_{\alpha_i} } ) $.
  Finally we define $ \pi( \lambda_{\alpha_0} )$ to be the operator $T(W_0(-1))$ defined in (3.7). 
\par It remains to prove that the operators defined above satisfies the relations (1)-(8) in Theorem 4.5. In the proof we will repeatedly use Lemma 3.2. 
\par The relation (1) is trivial from definition. For relation (2), the only case that needs care is the case when $\widetilde h=\widetilde h_{\alpha_0}(r)$ for some $r\in F^*$ since $\ol{h_{\alpha_i}(s)}\in P$ for $i\neq 0$. For $h=h_{\alpha_0}(r)$, notice that the matrix $\begin{pmatrix}
		r^{-1}&\\
		&r
	\end{pmatrix}$ appeared in formula (3.8) lies in $P$, so the relation is clear. 
\par For (3), first note that the two parts of (3) are equivalent because $\dot w_a$ normalizes $U^a$, so we only need to prove one. Also the only non-trivial case is $a=\alpha_0$ since $\ol{\lambda_{\alpha_i}}\in P$ for $i\neq0$. Let $u'=\dot w_{\alpha_0}u\dot w_{\alpha_0}^{-1}$, then we need to prove $$T(W_0(-1))T(\ol{\lambda(u)})=T(\ol{\lambda(u')})T(W_0(-1))$$ By the relations in the Kac-Moody group $G_D(F)$, we have $$W_0(-1)\ol{\lambda(u)}=\ol{\lambda(u')}W_0(-1)$$ Since $\ol{\lambda(u)},\ol{\lambda(u')}\in P$, by Lemma 3.2, we have $$T(\ol{\lambda(u')})T(W_0(-1))=T(\ol{\lambda(u')}W_0(-1))$$ if we take the same Haar measures on the images, and $$T(W_0(-1))T(\ol{\lambda(u)})=T(W_0(-1)\ol{\lambda(u)})$$ if we take the Haar measures on the images to be compatible with $\alpha_{\ol{\lambda(u)}}$. But since $u$ is unipotent, the induced map on finite-dimensional subspaces $\alpha_{\ol{\lambda(u)}}:\gamma_{W_0(-1)}\to\gamma_{W_0(-1)\ol{\lambda(u)}}$ has determinant $1$, thus preserves the Haar measures, so the choice of Haar measure on $\gamma_{W_0(-1)\ol{\lambda(u)}}$ when defining $T(W_0(-1)\ol{\lambda(u)})$ coincides with the choice of Haar measure on $\gamma_{\ol{\lambda(u')}W_0(-1)}=\gamma_{W_0(-1)\ol{\lambda(u)}}$ when defining $T(\ol{\lambda(u')}W_0(-1))$. This implies that $$T(\ol{\lambda(u')})T(W_0(-1))=T(W_0(-1))T(\ol{\lambda(u)})$$
\par For (4), also the only non-trivial case is the case that $a=\alpha_0$, so by (3.7) we have \begin{align*}
		&(T(W_0(-1))^2f)(\cdots,a_{-2},a_{-1}|a_0,a_1,\cdots)\\=&\gamma(-1)\int_F\psi(-a_0y)(T(W_0(-1))f)(\cdots,-a_2,-a_1|y,a_{-1},a_{-2},\cdots)dy\\=&\gamma(-1)^2\int_F\int_F\psi(-a_0y-yz)f(\cdots,-a_{-2},-a_{-1}|z,-a_1,-a_2,\cdots)dydz\\=&\gamma(-1)^2\int_F\delta(a_0+z)f(\cdots,-a_{-2},-a_{-1}|z,-a_1,-a_2,\cdots)dz\\=&\gamma(-1)^2f(\cdots,-a_{-2},-a_{-1}|-a_0,-a_1,\cdots)
	\end{align*} On the other hand, by formula (3.6), we have $$T(H_0(-1))=\gamma(-1)\gamma(1)^{-1}f(\cdots,-a_{-2},-a_{-1}|-a_0,-a_1,\cdots)$$ so relation (4) follows from the property of the Weil index that $\gamma(r)\gamma(-r)=1$.
	\par For relation (5) one also only need to check the case that $a=\alpha_0$. We can suppose $h=h_{\alpha_i}(r)$. If $i=2,3,\cdots,n$, we have ${\frak s}_0^{-1}\ol{h_{\alpha_i}(r)}=\ol{h_{\alpha_i}(r)}$ and $\ol{h_{\alpha_i}(r)}\in P$, so we can argue like in the proof of (3) since $\ol{h_{\alpha_i}(r)}$ does not change the measure on $\gamma_{W_0(-1)}$. The most complicated is the case that $i=1$, which involves a computation on $4\times 4$ matrices, but indeed the computation on the $(1,l+1)$-coordinates and the $(2,l+2)$-coordinates are independent, and the result on the $(2,l+2)$ coordinates are trivial, so it suffices to compute the $(1,l+1)$-coordinates, which reduces to the computation in the case $i=0$, which is done in the following. If $i=0$, note that by Lemma 3.2, we have $$T(W_0(-1))=T\matrixx{}{1}{1}{}T(D_1)=T(D_{-1})T\matrixx{}{-1}{-1}{}$$ so by (3.8) we have \begin{align*}
		&T(W_0(-1))^{-1}T(H_0(r))T(W_0(-1))\\=&|r|^\frac{1}{2}\gamma(r^{-1})\gamma(1)^{-1}T(D_{-1})T\matrixx{}{1}{1}{}T\matrixx{r^{-1}}{}{}{r}T\matrixx{}{1}{1}{}T(D_1)\\=&|r|^\frac{1}{2}\gamma(r^{-1})\gamma(1)^{-1}T(D_{-1})T\matrixx{r}{}{}{r^{-1}}T(D_1)
	\end{align*} Then by the definition of $T(D_1)$ and $T(D_{-1})$ in Section 3.2.1, we have \begin{align*}
		&(T(W_0(-1))^{-1}T(H_0(r))T(W_0(-1))f)(\cdots,a_{-2},a_{-1}|a_0,a_1,\cdots)\\=&|r|^\frac{1}{2}\gamma(r^{-1})\gamma(1)^{-1}(T(D_{-1})T\matrixx{r}{}{}{r^{-1}}T(D_1)f)(\cdots,a_{-2},a_{-1}|a_0,a_1,\cdots)\\=&|r|^\frac{1}{2}\gamma(r^{-1})\gamma(1)^{-1}\int_F(T\matrixx{r}{}{}{r^{-1}}T(D_1)f)(\cdots,-a_{-2}.-a_{-1},y_0|a_1,a_2,\cdots)\psi(y_0a_0)dy_0\\=&|r|^\frac{1}{2}\gamma(r^{-1})\gamma(1)^{-1}\int_F(T(D_1)f)(\cdots,-ra_{-2}.-ra_{-1},ry_0|r^{-1}a_1,r^{-1}a_2,\cdots)\psi(y_0a_0)dy_0\\=&|r|^\frac{1}{2}\gamma(r^{-1})\gamma(1)^{-1}\int_F\int_F f(\cdots,ra_{-2}.ra_{-1}|y_{-1},r^{-1}a_1,r^{-1}a_2,\cdots)\psi(y_0a_0)\psi(-ry_{-1}y_0)dy_0dy_{-1}\\=&|r|^{-\frac{1}{2}}\gamma(r^{-1})\gamma(1)^{-1}\int_F\int_F f(\cdots,ra_{-2}.ra_{-1}|r^{-1}y_{-1},r^{-1}a_1,r^{-1}a_2,\cdots)\psi(y_0a_0)\psi(-y_{-1}y_0)dy_0dy_{-1}\\=&|r|^{-\frac{1}{2}}\gamma(r^{-1})\gamma(1)^{-1}\int_Ff(\cdots,ra_{-2}.ra_{-1}|r^{-1}y_{-1},r^{-1}a_1,r^{-1}a_2,\cdots)\delta(a_0-y_{-1})dy_{-1}\\=&|r|^{-\frac{1}{2}}\gamma(r^{-1})\gamma(1)^{-1}(T\matrixx{r}{}{}{r^{-1}}f)(\cdots,a_{-2},a_{-1}|a_0,a_1,\cdots)\\=&(T(H_0(r^{-1}))f)(\cdots,a_{-2},a_{-1}|a_0,a_1,\cdots)
	\end{align*} Note that in the last step we used the fact that $\gamma(r)=\gamma(r^{-1})$. 
	\par For relation (6) one only need to check the Braid relations related to $\lambda_{\alpha_0}$. By the explicit description of the root datum in Section 2, the only such relation is that $\lambda_{\alpha_0}\lambda_{\alpha_1}\lambda_{\alpha_0}\lambda_{\alpha_1}=\lambda_{\alpha_1}\lambda_{\alpha_0}\lambda_{\alpha_1}\lambda_{\alpha_0}$. Let $w=\ol{\lambda_{\alpha_0}\lambda_{\alpha_1}\lambda_{\alpha_0}}$, fix a choice of Haar measures on $\gamma_{w}$, then it suffices to prove $T(w)T(\ol{\lambda_{\alpha_1}})=T(\ol{\lambda_{\alpha_1}})T(w)$ since $T(w)$ differs from $T(\ol{\lambda_{\alpha_0}})T(\ol{\lambda_{\alpha_1}})T(\ol{\lambda_{\alpha_0}})$ by a scalar. Note that $\ol{\lambda_{\alpha_1}}\in P$ is an element of finite order, thus it has determinant $1$ on any finite-dimensional subspace, so we can proceed as in the proof (3) to conclude that (6) is also correct. 
	\par For relations (7) and (8), still we only need to check the case $a=\alpha_0$, which follows from the process of defining the operators $T(W_0(r))$ and $T(H_0(r))$ in Section 3.2.

\

\


\begin{thebibliography}{CERP}


\bibitem[Bo]{Bo}N. Bourbaki, \textit{Groupes et Algebres de Lie}, Chaps IV-VI, Hermann, Paris, 1968. 


\bibitem[Gi]{Gi}  A. Givental, Gromov-Witten Invariants and Quantization of Quadratic Hamiltonians, {\it Moscow Math Journal}, 1(4), 2001.



\bibitem[GZ1]{GZ1}

H. Garland, Y. Zhu, {\it On the Siegel-Weil theorem for loop groups. I},
Duke Math. J. {\bf 157} (2011), no. 2, 283--336.

\bibitem[GZ2]{GZ2}

\bysame, {\it On the Siegel-Weil theorem for loop groups. II},
Amer. J. Math. {\bf 133} (2011), no. 6, 1663--1712.




\bibitem[K]{K} V. Kac, Infinite dimensional Lie algebras,  Cambridge University Press (1990).



\bibitem[LZ]{LZ}

D. Liu, Y. Zhu, {\it On the theta functional of Weil representations for symplectic loop groups},
J. Algebra {\bf 324} (2010), no. 11, 3115--3130.

\bibitem[Ma]{Ma} Hideya Matsumoto, Sur les sous-groupes arithmetiques des groupes semi-simples deployes, Ann. Sci. Ecole Norm. Sup. (4) 2
(1969), 1-62 (French). MR0240214





\bibitem[PP]{PP}
M. Patnaik, A. Puskas,  {\it Metaplectic Covers of Kac-Moody Groups and Whittaker functions,} Duke. J. Math. {\bf 168, } no. 4 (2019),  553-653.




\bibitem[Tits]{Tits}	J. Tits, Uniqueness and presentation of Kac-Moody groups over fields {\it J. Algebra }  105, no. 2 (1987): 542-73.


\bibitem[W]{W}

A. Weil, {\it Sur Certaines Groups d'operators unitaires}, Acta Math.
  {\bf 11} (1964), 143--211.



\bibitem[Z]{Z}

Y. Zhu, {\it Theta functions and Weil representations of loop symplectic groups}, Duke Math. J. {\bf 143} (2008), no. 1, 17--39.


\end{thebibliography}
\end{document}